\documentclass[review]{elsarticle}

\usepackage{lineno,hyperref}
\modulolinenumbers[5]

\journal{arXiv}
\usepackage{amscd}
\usepackage{mathrsfs}
\usepackage{amsfonts}
\usepackage{graphicx}
\usepackage{amsmath,amscd,stmaryrd,amssymb,array, amsfonts,amsmath,amssymb}
\usepackage{enumitem}
\usepackage{epstopdf}
\usepackage{caption}
\usepackage{subfigure}
\usepackage{float}

\numberwithin{figure}{section}
 \numberwithin{equation}{section}

\newtheorem{theorem}{Theorem}[section]
\newtheorem{proposition}[theorem]{Proposition}
\newtheorem{definition}[theorem]{Definition}

\newtheorem{lemma}[theorem]{Lemma}
\newtheorem{remark}[theorem]{Remark}

\allowdisplaybreaks[4]

\newcommand{\bi}{{\mathbf{i}}}

\newcommand{\bbC}{{\mathbb C}}
\newcommand{\bbE}{{\mathbb E}}
\newcommand{\bbG}{{\mathbb G}}

\newcommand{\bbM}{{\mathbb M}}
\newcommand{\bbN}{{\mathbb N}}

\newcommand{\cU}{{\mathcal U}}

\newcommand{\cS}{{\mathcal S}}

\newcommand{\mB}{\mb{B}}

\def\be{\begin{equation}}
\def\ee{\end{equation}}
\def\bes{\begin{equation*}}
\def\ees{\end{equation*}}
\def\bsp{\begin{split}}
\def\esp{\end{split}}

\def\ba{\begin{array}}
\def\ea{\end{array}}
\def\benu{\begin{enumerate}}
\def\eenu{\end{enumerate}}
\def\bt{\begin{theorem}}
\def\et{\end{theorem}}
\def\bp{\begin{proposition}}
\def\ep{\end{proposition}}
\def\bl{\begin{lemma}}
\def\el{\end{lemma}}
\def\br{\begin{remark}}
\def\er{\end{remark}}
\def\bd{\begin{definition}}
\def\ed{\end{definition}}

\def\b{\beta}
\def\De{\Delta}
\def\de{\delta}

\def\lam{\lambda}

\def\ve{\varepsilon}
\def\sig{\sigma}

\def\gam{\gamma}

\def\a{\alpha}

\def\.{\cdot}

\def\bbC{\mathbb{C}}

\def\R{\mathbb{R}}

\def\A{\forall}
\def\ol{\overline}

\def\ra{\rightarrow}
\def\stac{\stackrel}
\def\~{\tilde}
\def\8{\infty}
\def\X{\times}

\def\[{\left[}
\def\]{\right]}
\def\({\left(}
\def\){\right)}
\def\llg{\left\langle}
\def\rrg{\right\rangle}

\def\mb{\mbox}
\def\bmt{\begin{matrix}}
\def\emt{\end{matrix}}

\def\sm{\setminus}

\def\bx{\blacksquare}

\def\hs{\hspace{0.5cm}}
\def\Vs{\vskip8pt}\def\vs{\vskip4pt}

\def\({\left(}\def\){\right)}

\begin{document}

\begin{frontmatter}

\title{New Schemes for
 Solving the   Principal\\ Eigenvalue Problems of Perron-like  Matrices via Polynomial Approximations of Matrix Exponentials$^\dag$\footnote{$\dag$This work was supported by the National Natural Science Foundation of China [11871368]}}

\author[mymainaddress]{Desheng  Li\,\corref{mycorrespondingauthor}}
\cortext[mycorrespondingauthor]{Corresponding author}
\ead{lidsmath@tju.edu.cn}

\author[mymainaddress]{Ruijing Wang}
\ead{wrj-math@tju.edu.cn}
\address[mymainaddress]{School of Mathematics,  Tianjin University, Tianjin 300072,  China}

\begin{abstract}A  real square matrix is {Perron-like}  if it has a real eigenvalue $s$, called the {principal eigenvalue} of the matrix, and $\mb{Re}\,\mu<s$ for any other eigenvalue $\mu$.  Nonnegative matrices and symmetric ones are  typical examples of this class of matrices. The  main purpose of  this paper is to develop a set of new schemes to compute the principal eigenvalues of Perron-like  matrices   and the  associated generalized eigenspaces  by using  polynomial approximations of matrix exponentials.  Numerical examples show   that these schemes   are effective   in practice. \end{abstract}

\end{frontmatter}

\newpage
\tableofcontents
\newpage

\section{Introduction}

The study of numerical eigenvalue problems of matrices has witnessed a long history of more than one hundred years. By far numerous efficient   computation methods and  algorithms have been developed such as   the power method, the Jacobi method,  the LR and the   QR algorithms, the Lanczos procedure, and  the  Krylov subspace methods,
which  were systematically summarized and  discussed in many excellent books and monographs (see e.g. \cite{Bj,Ford,
GPS,GB,Hous,Lyche,Sad1,Wat1,Wat2,Wat3,Wilk}). 
These methods and algorithms enable us to acquire effectively  approximate eigenvalues and eigenvectors  of a given  matrix with prescribed  accuracy.  The subspace methods even allow  one  to compute invariant subspaces of matrices.
Of course, each method has   more or less some  limitation or drawbacks. For instance,
 the  convergence of the power method  is only guaranteed in a generic sense with respect to the initial vector $x_0$; furthermore, the choice of $x_0$ significantly affects
 the  convergence speed of the iteration sequence.   Similar situation also occurs   in the subspace methods, where  it is   often required  that the {\em initial} subspace $\cS$ in the iteration procedure satisfies
 \be\label{e:1.1}\cS\cap\, \cU^\perp=\{0\}\ee to guarantee the convergence, where $\cU^\perp$ denotes the orthogonal complement of the target subspace  $\cU$; see Watkins \cite[Theorems 5.1.1]{Wat2} and  Bj\"{o}rck \cite[Theorem 3.3.5]{Bj}. It is known that  \eqref{e:1.1} is true in  a generic sense with respect to $\cS$. However in practice, since  the target subspace $\cU$ is in fact unknown,   in many cases  we have   to try  our luck whether this requirement is fulfilled by a given subspace $\cS$. Let us also mention that most of the existing computation methods and algorithms  work well only for semisimple matrices whose  eigenvalues share the same algebraic and the geometric multiplicities.

In this present work, inspired by the dynamical approach towards the Perron-Frobenius theory  in \cite{LiDS},  we develop in a self-contained manner some  new schemes  for computing the principal  eigenvalues and the  corresponding generalized eigenspaces of a wide  class of real matrices which we call {\em Perron-like matrices} by using matrix exponentials and their polynomial approximations.  This  class of matrices contain nonnegative matrices and symmetric ones  as typical examples.
Similar as in the case of  the subspaces methods, our schemes make use of matrix iterations. One  advantage of the  schemes is that the convergence  is always guaranteed as long as the initial matrices are taken nonsingular. Another one  is that they can successfully generate the whole generalized principal eigenspace of a Perron-like matrix in  the non-semisimple  case. 

Now we give a more detailed description of our work and the organization of this paper.  In Section \ref{s:2} we make some preliminaries. In Section \ref{s:3} we prove some convergence results for Perron-like matrices.
 Denote by $\bbM_m$  the space consisting  of  $m\X m$ real matrices.
    A matrix $A\in \bbM_m$ is called {\em Perron-like} if it has a real eigenvalue $s$ (which will be called the {\em principal eigenvalue} of $A$), moreover,  $\mb{Re}\,\mu<s$ for all other eigenvalues $\mu$.
Let  $A\in\bbM_m$ be a Perron-like matrix  with principal eigenvalue  $s$.  Given a {\em nonsingular} matrix $V\in \bbM_m$ with column vectors $v_i$, $i\in \{1,2,\cdots,m\}:=J$, put $$X(t)=e^{tA}V.$$ Let $y_i=\Pi_1 v_i$, where $\Pi_1$ is the projection from $\R^m$ to the generalized eigenspace $\mb{GE}_s(A)$ associated with $s$. Denote by $Q_{y_i}(t)$ the vector-valued {\em characteristic polynomial} of $y_i$ as defined in \eqref{ecp} below. Let $Q_Y(t)$ be  the matrix with column vectors $Q_{y_i}(t)$ ($i\in J$). We show that  there exist  $B_0,\de>0$ such that
\be\label{e1.1}
\left\|\frac{X(t)}{\|X(t)\|}-\frac{Q_Y(t)}{\|Q_Y(t)\|}\right\|\leq B_0 e^{-\de t},\hs t\geq1.
\ee
Based on this fundamental result, we then verify that there is a matrix $\Xi$ whose column vectors $\xi_i\,\,(i\in J)$ span a nontrivial invariant  subspace of the eigenspace $\mb{E}_s(A)$ associated with $s$ such that
\be\label{e1.2}
\left\|\frac{X(t)}{\|X(t)\|}-\Xi\right\|\leq O(t^{-1})+B_1e^{-\de t},\hs t\geq1.
\ee
As a consequence, one naturally has
\be\label{e1.3}\left|\frac{1}{\|X(t)\|^2}\langle AX(t),X(t)\rangle-s\right|\leq O(t^{-1})+ B_2e^{-\de t},\hs t\geq 1.
\ee
If  $s$ is {\em semisimple}, i.e., $s$ shares the same algebraic and the geometric multiplicities (hence $\mb{GE}_s(A)=\mb{E}_s(A)$), then \be\label{e1.5}\mb{E}_s(A)=\mb{span}\{\xi_i:\,\,i\in J\}.\ee Furthermore,  the term  $O(t^{-1})$ in the above estimates can be removed, and therefore  the convergence in \eqref{e1.2} and \eqref{e1.3} is actually  exponential.

It is interesting to mention that  \eqref{e1.2} also leads to a strengthened version of the Perron theorem for nonnegative matrices; see Section \ref{s:3.3} for details.

Nowadays there are many effective  ways for computing    matrix exponentials;  see e.g. \cite{FH,HL,ML}. One of the simplest ways is  to use Taylor expansions which yield polynomial approximations of matrix exponentials. In Section \ref{s:4} we will   give some convergence  results parallel to those in \eqref{e1.1}-\eqref{e1.3}
for polynomial approximations of $X(t)$. Based on these results,  we then  design in Section \ref{s:5} a corresponding  iterative  scheme. Specifically, for each $n\in \bbN$, we define a mapping $K_n$ on $\Sigma_1:=\{X\in \bbM_m:\,\,\|X\|=1\}$ as
$$ K_n X=\frac{T_n X}{\|T_n X\| },\hs X\in \Sigma_1,$$
where $T_n=\sum_{\ell=0}^n\frac{1}{\ell!}A^\ell$.
Given a {nonsingular} matrix   $V\in \Sigma_1$, set
$$
X_n(0)=V,\hs\,X_n(k+1)=K_nX_n(k)\,\,\,\,(k,n\in\bbN).
$$
We show that there exist positive constants $B_0,B_1$ and $B_2$ such that
\be\label{e1.4}
\ba{ll}
\left\|X_n(n)-\frac{Q_Y(n)}{\|Q_Y(n)\|}\right\| \leq B_0e^{-\de n}+\(\frac{\lam\|A\|}{n+1}\)^{n+1}o(1),
\ea
\ee
\be\label{e:1.3}
\|X_n(n)-\Xi\| \leq O(n^{-1})+B_1e^{-\de n}+\(\frac{\lam\|A\|}{n+1}\)^{n+1}o(1),
\ee and
\be\label{e:1.4}\left|\langle AX_n(n),\,X_n(n)\rangle -s\right|\leq O(n^{-1})+B_2e^{-\de n}+ {\|A\|}\,\(\frac{\lam\|A\|}{n+1}\)^{n+1}o(1),
\ee
where $\lam=2e^{2\|A\|}$, and $o(1)$ is an infinitesimal as $n\ra\8$.

As in \eqref{e1.2} and \eqref{e1.3}, if  $s$ is semisimple then the first terms  in the righthand sides of  \eqref{e:1.3} and \eqref{e:1.4} can be dropped. In view of \eqref{e1.5} we see that   the above iterative scheme readily provides  an efficient way for solving the principal eigenvalue problems  of Perron-like matrices.
The situation in the non-semisimple case seems to be complicated. On one hand, \eqref{e:1.3} and \eqref{e:1.4} as well as  the numerical simulation in Section \ref{s:5} indicate that  the convergence can be very  slow, which fact may make the scheme to be of little practical sense. On the other hand, the column vectors of $\Xi$  only span a proper subspace of the generalized eigenspace $\mb{GE}_s(A)$ associated with $s$. 
In Sections \ref{s:6}\,-\,\ref{s:8} we develop some  new methods to deal with the  non-semisimple case.

In Section \ref{s:6} we design  a numerical scheme  to determine the {\em cyclic order} 
 of $\mb{GE}_s(A)$, namely, the smallest number $\nu$ such that
$$
(A-s I)^\nu \(\mb{GE}_s(A)\)=\{0\}.
$$
Because the  scheme  is mainly  based on the exponential convergence in \eqref{e1.4} and a ``\,rough'' approximate value $\bar s$ of $s$, it  can help us determine the cyclic order of $\mb{GE}_s(A)$ quickly. This allows us to develop a more efficient  combined method to compute the principal eigenvalue $s$, which will be addressed in Section \ref{s:7}.

Section \ref{s:8} is concerned with  the computation of the
   generalized eigenspace $\mb{GE}_s(A)$.
   Let    $V\in \bbM_m$ be a nonsingular matrix with column vectors  $v_i\,\,({i\in J})$.
 Then one  trivially verifies that   $\mb{GE}_s(A)=\mb{span}\{y_i:\,\,i\in J\}$, where $y_i=\Pi_1 v_i$.
  Therefore the problem of computing  $\mb{GE}_s(A)$ reduces to that  of the matrix $Y$  with column vectors $y_i$ $(i\in J)$. Based on the simple observation that
$$\ba{ll}
Y
=P(t)Q_{Y}(t)=P(t)\hat X(t)-P(t)\hat Z(t),\ea
$$
where $\hat X(t)=e^{t(A-sI)}V$,\, $\hat Z(t)=e^{t(A-sI)}(V-Y)$, and
$$
P(t)=\sum_{k=0}^d(-1)^k\,\,\frac{t^k}{k!}(A-sI)^k, \hs \nu-1\leq d\leq m-1,$$ in Section \ref{s:8} we propose  an iterative scheme for computing the matrix $Y$ by using polynomial approximations of matrix exponentials.

Numerical examples will be given in Sections \ref{s:5}\,-\,\ref{s:8} to illustrate  the efficiency of the  computation schemes.

\section{Preliminaries}\label{s:2}
This  section is concerned with some preliminaries.  
\subsection{Basic notions and notations}
\vs \noindent$\bullet$ {\bf Notations} \,\,Throughout the paper $\R$ (resp. $\bbC$) stands for the field of real (resp. complex) numbers. Let  $\R_+=[0,\8)$, and $\bbN$ the set of nonnegative integers.
       Denote by $\R^m$  the Euclid space consisting of  $m$-dimensional  {\em column vectors}  equipped with the usual inner product $\langle\.,\,\.\rangle$ and the  $2$-norm $\|\.\|$.
Let $\bbM_m:=\R^{m\X m}$ be the space of  $m\X m$ real matrices.
 The same notations $\langle\.,\,\.\rangle$ and $\|\.\|$ as above will be used to denote the    inner product and the $2$-norm of $\bbM_m$, respectively.
A matrix $A\in\bbM_m$ with  column vectors $\a_i$ ($i\in J$) will be written  as $A=[\a_i]_{i\in J}$, where (and below) $J$ denotes the index set $\{1,2,\cdots,m\}$. It is obvious   that
$$
\langle A,B\rangle=\sum_{i\in J}\langle\a_{i},\,\b_{i}\rangle,\hs\,\A\,A=[\a_i]_{i\in J},\,B=[\b_i]_{i\in J}\in\bbM_m.
$$
Denote by $I$ the identity matrix in $\bbM_m$.

\vs
\noindent$\bullet$ {\bf Spectrum and invariant subspaces of matrices} \,\,Given  $A\in\bbM_m$, denote by   $\sig(A)$    the {\em spectrum} of $A$ consisting of the eigenvalues of   $A$.
    The {\em spectral bound}  $s(A)$ and {\em spectral radius} $r(A)$ of $A$ are defined  as
$$
s(A)=\max\{\mb{Re}\,\mu:\,\,\mu\in\sig(A)\},\hs r(A)=\max\{|\mu|:\,\,\mu\in\sig(A)\}.
$$
One can naturally think of $A$ as a linear transform on either $\R^m$ or $\bbC^m$. A subspace $Y$ of $\R^m$ (or $\bbC^m$) is called an {\em invariant subspace} of $A$, if $AY\subset Y$.

 For each   $\mu\in\sig( A)$, the space
\be\label{e:es0}
\bbG\bbE_\mu(A):=\{\xi\in\bbC^m:\,\,( A-\mu I)^{k}\xi=0\mb{ for some }k\geq 1\}
\ee
is  a nontrivial  invariant subspace of $A$.
Consequently
\be\label{e:2.5}\mb{GE}_\mu( A):=\{x,\,y:\,\,u=x+\bi y\in \bbG\bbE_\mu(A)\}\ee is an  invariant subspace of $A$ (in $\R^m$), where $\bi$ stands for the unit imaginary number in $\bbC$. We call
$\mb{GE}_\mu( A)$    the (real) {\em generalized  eigenspace} (or, {\em singular space}) of $A$ associated with $\mu$.

\Vs
\noindent$\bullet$ {\bf The case of  real eigenvalues} \,\,Now we assume that $\mu$ is a real eigenvalue of
$A\in \bbM_m$. Then  $( A-\mu I)^{k}(x+\bi y)=0$ amounts to say that $( A-\mu I)^{k}x=0=( A-\mu I)^{k}y$. It follows by \eqref{e:es0} and \eqref{e:2.5} that
 \be\label{ss}
\mb{GE}_\mu( A)=\{\xi\in\R^m:\,\,( A-\mu I)^{k}\xi=0\mb{ for some }k\geq 1\}.
\ee
We call each nonzero element $\xi$ in $\mb{GE}_\mu( A)$ a {\em generalized eigenvector} of $A$. Let $\xi$ be a generalized eigenvector of $A$ corresponding to  $\mu$. Then there is an integer $\nu\geq 1$ such that
\be\label{e:2.1}
(A-\mu I)^j\xi\ne0\,\,\,(0\leq j\leq \nu-1),\hs (A-\mu I)^\nu\xi=0;\ee
furthermore,  $w:=(A-\mu I)^{\nu-1}\xi$ is an eigenvector of $A$. For convenience, we call $\nu$ the {\em order} \,of $\xi$, denoted  by $\mb{Ord}\,(\xi)$. Define   $$\ba{ll}\mb{Cord}\(\mb{GE}_\mu(A)\):=\max_{\xi\in \mb{\footnotesize GE}_\mu(A)}\mb{Ord}(\xi).\ea$$  $\mb{Cord}\(\mb{GE}_\mu(A)\)$ is called the {\em cyclic order} of $\mb{GE}_\mu(A)$.

Denote by $\mb{E}_\mu(A)$ the {\em eigenspace} of $A$  spanned by all the eigenvectors of $A$ corresponding  to $\mu$. An  eigenvector $w\in \mb{E}_{\mu}(A)$ is called a {\em dominant eigenvector}, if there is  a $\xi\in\mb{GE}_{\mu}(A)$ such that $$w=(A-\mu I)^{\nu-1}\xi,\hs \mb{where } \nu=\mb{Cord}\,\(\mb{GE}_{\mu}(A)\).$$
  The {\em dominant eigenspace}  of $A$ associated with $\mu$, denoted by  $\mb{DE}_\mu(A)$, is the subspace of $\mb{E}_\mu(A)$ spanned by  the dominant   eigenvectors of $A$ corresponding to $\mu$.
Clearly
$$
\mb{DE}_\mu(A)\subset \mb{E}_\mu(A)\subset \mb{GE}_\mu(A).
$$

The validity of the following simple fact is almost obvious.
\bl\label{l:2.3}$
(A-\mu I)^{\nu-1}\mb{\em GE}_\mu(A)=\mb{\em DE}_\mu (A),
$
 $\mb{where } \nu=\mb{\em Cord}\,\(\mb{\em GE}_{\mu}(A)\).$
\el
\bd
We say that  $\mu$ is  {semisimple }, if $\mb{\em GE}_\mu(A)=\mb{\em E}_\mu(A)$.
\ed
It is easy to see that if $\mu$ is semisimple, then
\be\label{e:2sm}
\mb{DE}_\mu(A)=\mb{E}_\mu(A)= \mb{GE}_\mu(A).
\ee

\vs \noindent$\bullet$ {\bf The characteristic polynomial of a generalized eigenvector} \,\,
\vs
Let
$A\in \bbM_m$, and let  $\mu$ be a real eigenvalue of $A$. Given $\xi\in \mb{GE}_{\mu}(A)$, set \be\label{ecp}Q_\xi(t):=\sum_{k=0}^{\nu_\xi-1}\frac{t^k}{k!}(A-\mu I)^k\xi,\hs\mb{where }\,\nu_\xi=\mb{Ord}(\xi).\ee $Q_\xi(t)$ is  called the {\em characteristic   polynomial} of $\xi$. Since $(A-\mu I)^k\xi=0$ for $k\geq \nu_\xi$, we obviously have
  \be\label{ecp2}Q_\xi(t)=\sum_{k=0}^{d}\frac{t^k}{k!}(A-\mu I)^k\xi=e^{t(A-\mu I)}\xi\ee
  for any integer $d\geq \nu_\xi-1$.
 \br\label{r:2.1}
   It is clear that if $\xi\in \mb{\em E}_\mu(A)$ then $Q_\xi(t)\equiv \xi$. \er
\subsection{Perron-like matrices}
\vs

 \bd A matrix  $A\in \bbM_m$ is said to be   Perron-like, if
the spectral bound  $s:=s(A)\in\sig(A)$; furthermore, $\mb{\em Re}\,\mu<s$ for all $\mu\in\sig(A)\sm\{s\}$.\ed

 For a Perron-like matrix $A$, we call  $s:=s(A)$ the {\em principal} eigenvalue of $A$ with  the invariant subspaces  $ \mb{ GE}_{s}(A)$, $\mb{E}_s(A)$ and $\mb{DE}_s(A)$ being referred to as   the   {\em principal} generalized eigenspace,  {eigenspace},  and dominant eigenspace of $A$, respectively.
Correspondingly,  each nonzero element in $ \mb{ GE}_{s}(A)$ (resp. $\mb{E}_s(A)$, $\mb{DE}_s(A)$) is called   a {\em principal} generalized eigenvector (resp. { eigenvector}, dominant eigenvector).

Now let us   give several examples of Perron-like matrices.
\Vs
\noindent{\bf Example 2.1.} A symmetric matrix $A\in \bbM_m$ is  Perron-like  with each eigenvalue being semisimple .
\Vs
\noindent{\bf Example 2.2.} Every  nonnegative matrix $A\in \bbM_m$ is a Perron-like one with  $s(A)=r(A)$.
This follows from some classical Krein-Rutman type theorems for matrices and bounded linear  operators; see e.g. \cite[Theorem 4.1]{LiDS}.
\Vs
\noindent{\bf Example 2.3.}  If all the off-diagonal entries of a matrix  $A\in \bbM_m$  are nonnegative, then $A$ is a Perron-like matrix. Indeed, in such a case one can pick a positive number $a>0$ sufficiently large so that $aI+A$ is nonnegative, from which and Example 2.2 one immediately concludes that $A$ is a Perron-like matrix with principal eigenvalue $$s:=s(A)=r(aI+A)-a.$$

\br One can easily give examples of Perron-like matrices that  are not covered by Examples 2.1-2.3.\er

\subsection{On the growth rate of  exponential functions of matrices}
\vs
The following lemma contains  a  well known fundamental fact   concerning the growth rate of exponential functions of matrices,  which actually holds true for general  bounded linear  operators and some types  of unbounded  operators in abstract  Banach spaces; see e.g. Henry \cite[Theorem 1.5.3]{Henry}. Here we include a proof for the reader's convenience.
\bl\label{l:2.1}
Let $A\in \bbM_m$. Assume that $s(A)\leq \lam$. Then for any $\ve>0$, there is a constant $C>0$ such that
\be\label{e:2.0}
\|e^{tA}\|\leq C e^{(\lam+\ve)t},\hs t\geq 0.
\ee
\el
{\bf Proof.}
We think of  $A$ as a linear transform (operator) on $\bbC^m$.
 Let $\mu=\a+\bi\b\in\sig(A)$, and  $\xi\in \bbG\bbE_\mu(A)$.
 Then $( A-\mu I)^\nu\xi=0$ for some $\nu\geq 1$. Hence
\be\label{e:2.6}
\ba{ll}e^{tA}\xi&=e^{\mu t}e^{t( A-\mu I) }\xi=e^{\mu t} \sum_{k=0}^\8 \frac{t^k}{k!}( A-\mu I)^k \xi\\[2ex]
&=e^{\mu t}Q_\xi(t)=e^{(\a+\ve)t}\left(e^{-\ve t}e^{\bi\b t}Q_\xi(t)\right),
\ea
\ee
where $Q_\xi(t)=\sum_{k=0}^{\nu-1} \frac{t^k}{k!}( A-\mu I)^k\xi.$   Since $Q_\xi(t)$ has only a polynomial growth, one easily sees  that there is $C_0>0$ such that
$
\left\|e^{-\ve t}e^{\bi\b t}Q_\xi(t)\right\|\leq C_0\|\xi\|
$
for $t\geq 0$. Thus by \eqref{e:2.6} we have
\be\label{e:2.2}
\|e^{tA}\xi\|\leq C_0e^{(\a+\ve)t}\|\xi\|\leq C_0e^{(\lam+\ve)t}\|\xi\|,\hs t\geq 0.
\ee

We infer from the basic knowledge in linear algebra that the space $E:=\bbC^m$ has a basis $\{\xi_1,\xi_2,\cdots,\xi_m\}$, where each $\xi_i$ belongs to some invariant subspace  $\bbG\bbE_{\mu_i}(A)$. For each $\xi_i$, by \eqref{e:2.2} there is a constant $C_i>0$ such that
\be\label{e:2.3}
\|e^{tA}\xi_i\|\leq C_ie^{(\lam+\ve)t}\|\xi_i\|,\hs t\geq 0.
\ee
Now we  write each $u\in E$ as $u=\sum_{1\leq i\leq m}a_i\xi_i$, where $a_i\in\bbC$. Define $|||u|||=\sum_{1\leq i\leq m}|a_i|$. Then $|||\.|||$ is an equivalent norm on $E$. Let $C'=\max_{1\leq i\leq m}C_i$.
Thanks to \eqref{e:2.3}, we deduce that
\be\label{e:2.4}\ba{ll}
\|e^{tA}u\|&=\left\|\sum_{1\leq i\leq m}a_ie^{tA}\xi_i\right\|\leq C'e^{(\lam+\ve)t}|||u|||\\[2ex]
&\leq C''e^{(\lam+\ve)t}\|u\|,\hs \A\,u\in E,\,\,t\geq 0,
\ea\ee
from which  we see that the {\em operator norm} of $e^{tA}$ (as an operator on $\bbC^m$) is dominated by $C''e^{(\lam+\ve)t}$. Consequently  the operator norm of $e^{tA}$ (as an operator on $\R^m$) is dominated by $C''e^{(\lam+\ve)t}$ as well. The estimate in \eqref{e:2.0} then immediately  follows because all the norms on a finite-dimensional Banach space  are equivalent. $\bx$
\Vs
\section{Fundamental  Convergence Results for Perron-like Matrices}\label{s:3}
\vs
In this section we prove some fundamental  convergence  results  mentioned in the introduction for Perron-like matrices, which serve  as the starting point of this work.
As an  interesting corollary, we also give a generalized Perron-like theorem for nonnegative matrices.

 The following notations will be used throughout the paper.
\benu\item[$\bullet$] The notation   $O(\ve)$ will stand for   a general infinitesimal  as $\ve\ra 0$ that can be  dominated by $C|\ve|$ for some $C>0$.
 \item[$\bullet$]   We denote by  $X^{\mb{\tiny T}}$ the transpose of a vector (or, matrix) $X$.
\eenu
\subsection{Preliminaries}\label{s:3.0}

\vs
 For simplicity, we  write  $E:=\R^m$, and let $J=\{1,2,\cdots,m\}$.
Assume that $A\in \bbM_m$ is  a Perron-like matrix with the principal eigenvalue $s:=s(A)$.

 Set $\sig_1=\{s\}$, and $\sig_2=\sig(A)\sm\{s\}$.
Then  $E$ has a direct sum decomposition   $E=E_1\oplus E_2$ corresponding to the spectral decomposition   $\sig(A)=\sig_1\cup \sig_2$,
\be\label{e:ds}\ba{ll}E_1=\mb{GE}_s(A),\hs E_2=\bigoplus_{\mu\in\sig_2}\mb{GE}_\mu(A).\ea\ee
Denote by $\Pi_j$ ($j=1,2$) the projection from $E$ to $E_j$.

\vs Let $\{v_i\}_{i\in J}$ be a basis  of $E$.
First, we have the following easy lemma.

\bl\label{l:3.1}$
\mb{\em span}\{\Pi_j v_1,\cdots,\Pi_j v_m\}=E_j,$ $j=1,2.$
\el
{\bf Proof.} Suppose, say, that $\mb{span}\{\Pi_1v_1,\cdots,\Pi_1v_m\}:=E_0\varsubsetneq E_1$. Then
$$\ba{ll}\mb{span}\{v_1,\cdots,v_m\}&\subset \bigoplus_{j=1,2}\mb{ span}\{\Pi_j v_1,\cdots,\Pi_j v_m\}\\[2ex]
&\subset E_0\oplus E_2\ne E,\ea
$$
which leads to a contradiction. $\bx$
\Vs

Let $$V=[v_i]_{i\in J},\hs Y=[y_i]_{i\in J},\hs Z=[z_i]_{i\in J},$$ where
$$y_i=\Pi_1v_i,\hs z_i=\Pi_2v_i.$$   Set $$
Q_Y(t)=\[Q_{y_i}(t)\]_{i\in J},\hs t\geq 0,
$$
where $Q_{y_i}(t)$ is  the (vector-valued) characteristic polynomial of $y_i$ as defined in \eqref{ecp2}. Put
\be\label{e:3.0}\mb{$x_i(t)=e^{tA}v_i$\, ($i\in J$),}\hs X(t)=\[x_i(t)\]_{i\in J}.\ee
It is trivial to see  that $X(t)=e^{tA}V$.

\subsection{Some fundamental  convergence theorems}\label{s:3.1}
\vs
We  now state and prove the  first convergence result.

\bt\label{t:3.1} There exist $B_0,\de>0$
such that
\be\label{e3.1}
\left\|\frac{x_i(t)}{\|X(t)\|}-\frac{Q_{y_i}(t)}{\left\|Q_Y(t)\right\|}\right\|\leq B_0 e^{-\de t},\hs t\geq1,\,\,i\in J.
\ee
\et

\noindent{\bf Proof.} Let $A_2=A|_{E_2}$ be the restriction of $A$ on $E_2$. Clearly  $\sig(A_2)=\sig_2$. Noticing that   $\mb{Re}\,\mu<s$ for all $\mu\in \sig(A_2)$, one can pick a $\de>0$ such that
\be\label{e:3.1}\mb{Re}\,\mu\leq s-2\de,\hs \mu\in \sig(A_2).\ee Lemma \ref{l:2.1} then asserts that  there is $C_1>0$ such that
\be\label{e3.7}
\|e^{tA_2}\|\leq C_1 e^{(s-\de)t},\hs t\geq 0.
\ee

For each  $i\in J$, we observe that
\be\label{e:5.10}\ba{ll}x_i(t)&=e^{tA}v_i=e^{tA}y_i+e^{tA}z_i\\[2ex]
&=e^{st}\sum_{k=0}^\8 \frac{t^k}{k!}( A-s I)^ky_i+e^{tA_2}z_i\\[2ex]
&=e^{st}Q_{y_i}(t)+ e^{tA_2}z_i,\hs t\geq 0.
\ea\ee
For notational simplicity, let us  write
$$\nu_i=\mb{Ord}\,(y_i),,\hs\mb{and }\,\nu=\mb{Cord}\(\mb{GE}_s(A)\).$$
Set $\~x_i(t)=e^{-st}t^{-(\nu-1)}x_i(t)$.
Then by \eqref{e:5.10} we have
\be\label{e3.3}
\~x_i(t)=p_i(t)+\~z_i(t),\hs t>0.
\ee
where
 \be\label{e:3.pi}p_i(t)=t^{-(\nu-1)}Q_{y_i}(t),\hs \~z_i(t)=e^{-st}t^{-(\nu-1)}e^{tA_2}z_i.\ee
 By \eqref{e3.7} one  finds that
\be\label{e3.8}\ba{ll}
\|\~z_i(t)\|&\leq C_1e^{-\de t}t^{-(\nu-1)}\|z_i\|\\[2ex]
&\leq C_1e^{-\de t}\|\Pi_2v_i\|\leq C_2\|v_i\|e^{-\de t},\hs t\geq 1,
\ea
\ee
where $C_2=C_1\|\Pi_2\|$. Thus if we  write
$$\|\~x_i(t)\|=\|p_i(t)\|+r_i(t),$$ then  by \eqref{e3.3} it is easy to deduce that
\be\label{e3.11}
|r_i(t)|\leq \|\~z_i(t)\|\leq C_2\|v_i\|e^{-\de t},\hs i\in J.
\ee

Let $\|\~X(t)\|:=e^{-s t}t^{-(\nu-1)}\|X(t)\|$ \,($t>0$). Obviously
\be\label{e:3.9a}
\|\~X(t)\|=\(\sum_{i\in J}\|\~x_i(t)\|^2\)^{1/2}=\(\sum_{i\in J}\(\|p_i(t)\|+r_i(t)\)^2\,\)^{1/2}.
\ee
It follows by the trigonal inequality of the  norm $\|\.\|$  that
\be\label{e:3.9b}\ba{ll}
\varrho(t)-r(t)\leq \|\~X(t)\|\leq \varrho(t)+r(t),\hs t\geq 1,
\ea\ee
where
\be\label{e:3.9c}\varrho(t)=\(\sum_{i\in J}\|p_i(t)\|^2\,\)^{1/2},\hs r(t)=\(\sum_{i\in J}|r_i(t)|^2\,\)^{1/2}.
\ee

We now evaluate  $\left\|\frac{x_i(t)}{\|X(t)\|}-\frac{p_i(t)}{\varrho(t)}\right\|$. For this purpose, let us write \be\label{e:3.18}\|\~X(t)\|=\varrho(t)+R(t).\ee Then by \eqref{e3.11} and \eqref{e:3.9b} it is easy to deduce that
\be\label{e3.29}\ba{ll}|R(t)|\leq r(t)\leq C_2\|V\| e^{-\de t},\hs t\geq 1.\ea\ee
 Observe that
\be\label{e:3.7}\displaystyle{\ba{ll}
\left\|\frac{x_i(t)}{\|X(t)\|}-\frac{p_i(t)}{\varrho(t)}\right\|&=\left\|\frac{\~x_i(t)}{\|\~X(t)\|}-\frac{p_i(t)}{\varrho(t)}\right\|=\left\|\frac{p_i(t)+\~z_i(t)}{\varrho(t)+R(t)}-\frac{p_i(t)}{\varrho(t)}\right\|\\[2ex]
&=\frac{\|\varrho(t)\~z_i(t)-R(t)p_i(t)\|}{\(\varrho(t)+R(t)\)\varrho(t)}.
\ea}
\ee
Since $\|p_i(t)\|\leq \varrho(t)$, we have
$$
 \|\varrho(t)\~z_i(t)-R(t)p_i(t)\|\leq \varrho(t)\(\|\~z_i(t)\|+|R(t)|\).
 $$
Therefore by \eqref{e:3.7} one gets that
\be\label{e3.13}\ba{ll}
\left\|\frac{x_i(t)}{\|X(t)\|}-\frac{p_i(t)}{\varrho(t)}\right\|&\leq \frac{\|\~z_i(t)\|+|R(t)|}{\varrho(t)+R(t)}=\frac{\|\~z_i(t)\|+|R(t)|}{\|\~X(t)\|}\\[2ex]
& \leq {(m+1)C_2}\|V\|\,\frac{1}{\|\~X(t)\|}e^{-\de t},\hs t\geq 1.
\ea
\ee

Set $J_s=\{i\in J:\,\,\nu_i=\nu\}$, and let
\be\label{e:3.14}
w_i=\frac{1}{(\nu-1)!}(A-sI)^{\nu-1}y_i,\hs i\in J_s.
\ee
We claim  that
\be\label{e:3.15}\mb{DE}_s (A)=\mb{span}\{w_i:\,\,\,i\in J_s\}.\ee
Indeed, by Lemma \ref{l:2.3} and Lemma \ref{l:3.1} we have
\be\label{e:3.16}\ba{ll}
\mb{DE}_s (A)&=(A-sI)^{\nu-1}\mb{GE}_s(A)=(A-sI)^{\nu-1}E_1\\[2ex]
&=(A-sI)^{\nu-1}\(\mb{span}\{y_i:\,\,\,i\in J\}\)\\[2ex]
&=\mb{span}\{(A-sI)^{\nu-1}y_i:\,\,\,i\in J\}\\[2ex]
&=\mb{span}\{(A-sI)^{\nu-1}y_i:\,\,\,i\in J_s\}\\[2ex]
&=\mb{span}\{w_i:\,\,\,i\in J_s\}.
\ea
\ee
Hence the claim holds true.

For convenience, we assign
\be\label{ew2}
w_i=0,\hs i\in J\sm J_s.
\ee
Then by  the definition of $Q_{y_i}(t)$ one can rewrite $p_i(t):=t^{-(\nu-1)}Q_{y_i}(t)$ as
\be\label{e:3.10}
p_i(t)=w_i+h_i(t),\hs i\in J.
\ee
where
$$
h_i(t)=\left\{\ba{ll}t^{-(\nu-1)}\sum_{k=0}^{\nu-2}\frac{t^k}{k!}(A-sI)^{k}y_i,\hs& i\in J_s;\\[2ex]
t^{-(\nu-1)}\sum_{k=0}^{\nu_i-1}\frac{t^k}{k!}(A-sI)^{k}y_i,&\mb{otherwise},\ea\right.
$$
In any case it can be easily seen that
\be\label{e:3.11}
\|h_i(t)\|\leq \|y_i\|O(t^{-1})\leq \|v_i\|O(t^{-1}).
\ee
Therefore by \eqref{e:3.18} and \eqref{e3.29} we deduce that
\be\label{e3.19}\lim_{t\ra\8}\|\~X(t)\|=\lim_{t\ra\8}\varrho(t)=\omega>(\mb{by \eqref{e:3.16}})>0,\ee
where
\be\label{dom}\ba{ll} \omega=\(\sum_{i\in J}\|w_i\|^2\)^{1/2}= \(\sum_{i\in J_s}\|w_i\|^2\)^{1/2}.\ea\ee
 Because $\|\~X(t)\|,\,\varrho(t)>0$ for all $t\geq 1$, \eqref{e3.19} implies that there is a positive number $\kappa>0$ such that
\be\label{e3.28}
\|\~X(t)\|, \,\varrho(t)\geq \kappa>0,\hs\,t\geq 1.
\ee
 Thus by \eqref{e3.13} one concludes  that
\be\label{e3.14}
\left\|\frac{x_i(t)}{\|X(t)\|}-\frac{p_i(t)}{\varrho(t)}\right\|\leq C_3\|V\| e^{-\de t},\hs t\geq 1,
\ee
where $C_3={(m+1)C_2}/{\kappa}$.
By \eqref{e:3.pi} and \eqref{e:3.9c} it can be easily seen that
\be\label{eq_i}\frac{p_i(t)}{\varrho(t)}=\frac{Q_{y_i}(t)}{\|Q_{Y}(t)\|},\ee
from which and   \eqref{e3.14} the convergence result in  \eqref{e3.1} immediately follows. $\bx$

\benu\item[$\bullet$]\,{\bf Notation} \,\,In the remaining part of this paper  the notation  $\de$  will always stand for the positive constant  given in \eqref{e:3.1} for a given Perron-like matrix $A$,
 which is determined by the spectral gap between the principal eigenvalue $s$ and the other  eigenvalues of $A$.
\eenu

Let us  proceed with the argument in the proof of Theorem \ref{t:3.1}.

Note that  $\varrho(t)$ (see \eqref{e:3.9c}) can be rewritten as  $\varrho(t)=\omega+\varrho_0(t),$ where  $\omega$ is number given in \eqref{dom}.   By \eqref{e:3.10}   and \eqref{e:3.11} one easily checks  that
 \be\label{e:3.6}
 |\varrho_0(t)|\leq\(\sum_{i\in J}\|h_i(t)\|^2\)^{1/2}\leq \|V\|O(t^{-1}).
 \ee

Now we  put
\be\label{e:3.3}
\xi_i=\omega^{-1}{w_i},\hs i\in J.
\ee
We show   that for each $i\in J$,
\be\label{e3.15}\left\|\frac{p_i(t)}{\varrho(t)}-\xi_i\right\|=O(t^{-1})\hs \mb{as }\,t\ra \8.\ee
 Indeed, by \eqref{e:3.11} and \eqref{e:3.6} we have
$$\ba{ll}\left\|\frac{p_i(t)}{\varrho(t)}-\xi_i\right\|
&=\left\|\frac{w_i+h_i(t)}{\omega+\varrho_0(t)}-\frac{w_i}{\omega}\right\|=\left\|\frac{\omega h_i(t)-\varrho_0(t)w_i}{\omega(\omega+\varrho_0(t))} \right\|\\[2ex]
&\leq \frac{\omega\|h_i(t)\|+|\varrho_0(t)|\|w_i\|}{\omega(\omega+\varrho_0(t))}\\[2ex]
&\leq \frac{\|h_i(t)\|+|\varrho_0(t)|}{\omega+\varrho_0(t)}=\frac{\|h_i(t)\|+|\varrho_0(t)|}{\varrho(t)} \\[2ex]
&\leq \kappa^{-1}\|V\|O(t^{-1})\hs(\mb{as $t\ra \8$}),\ea$$
where $\kappa$ is the constant given in \eqref{e3.28}. This is precisely what we want.

\vs Combining \eqref{e3.14} and  \eqref{e3.15}, it yields
\be\label{e3.16}\left\|\frac{x_i(t)}{\|X(t)\|}-\xi_i\right\|\leq O(t^{-1})+B_1e^{-\de t},\hs t\geq 1.\ee
We infer from  \eqref{e:3.15} that
   \be\label{e:3.25}\mb{DE}_s (A)=\mb{span}\{\xi_i:\,\,\,i\in J_s\}.\ee
It is also clear that
\be\label{e:3.26}
\sum_{i\in J}\|\xi_i\|^2=\sum_{i\in J_s}\|\xi_i\|^2=1.
\ee

Write $$ \psi_i(t)=\frac{x_i(t)}{\|X(t)\|},\hs i\in J.$$
Since  $\xi_i\in\mb{DE}_s(A)$ for all $i\in J$ (note that $\xi_i=0$ for $i\in J\sm J_s$), we have
$$\ba{ll}
&\,\,\,\,\,\left|{\sum_{i\in J}\langle A \psi_i(t),\, \psi_i(t)\rangle}-s\right|\\[2ex]
&=\left|{\sum_{i\in J}\langle A \psi_i(t),\, \psi_i(t) \rangle}-s\sum_{i\in J}\|\xi_i\|^2\right|\\[2ex]
&=\left|{\sum_{i\in J}\langle A \psi_i(t),\, \psi_i(t) \rangle}-\sum_{i\in J}\langle A\xi_i,\,\xi_i \rangle\right|\\[2ex]
&\leq \sum_{i\in J}\left|\langle A( \psi_i(t)-\xi_i),\, \psi_i(t)\rangle\right|+\sum_{i\in J}\left|\langle A\xi_i,\, \psi_i(t)-\xi_i \rangle\right|\\[2ex]
&\leq \sum_{i\in J}\|A( \psi_i(t)-\xi_i)\|+\sum_{i\in J}\|A\xi_i\|\,\| \psi_i(t)-\xi_i\|\\[2ex]
&\leq 2\|A\|\sum_{i\in J}\| \psi_i(t)-\xi_i\| \\[2ex]
&\leq(\mb{by \eqref{e3.16}})\leq  O(t^{-1})+2m B_1\|A\|e^{-\de t},\hs t\geq 1.
\ea$$

Summarizing  the above results, one obtains the following theorem.

\bt\label{t:3.3} There exist $B_1,B_2>0$ and $\xi_i\in \mb{\em E}_s(A)$ $(i\in J)$ with
\be\label{e:c3.20}\sum_{i\in J}\|\xi_i\|^2=1\,\,\mb{ and }\,\,\,\mb{\em DE}_s (A)=\mb{\em span}\{\xi_i:\,\,\,i\in J\}\ee
such that
\be\label{e:c3.1}
\left\| \frac{x_i(t)}{\|X(t)\|}-\xi_i\right\|\leq O(t^{-1})+B_1e^{-\de t},\hs t\geq1,
\ee
and
\be\label{e:c3.2}\left|\frac{1}{\|X(t)\|^2}\,{\sum_{i\in J}\langle A x_i(t),\, x_i(t)\rangle}-s\right|\leq O(t^{-1})+B_2e^{-\de t},\hs t\geq1.
\ee
\et

\subsection{The case of a semisimple  principal eigenvalue}\label{s:3.2}
\vs
 Now we consider the  particular case where  $ s$  is semisimple.
In such a case,  every  $y_i:=\Pi_1 v_i$ with $y_i\ne 0$ in the above argument is an eigenvector of $A$. Hence   $Q_{y_i}(t)\equiv y_i$ for $t\geq 0$. Consequently  \eqref{e3.1} reads as
 \be\label{e:3.12}
\left\|\frac{x_i(t)}{\|X(t)\|}-\xi_i\right\|\leq B e^{-\de t},\hs t\geq1,
\ee
where  $$\xi_i=\frac{y_i}{\|Y\|}.$$  This also leads to  a corresponding   strengthened version of \eqref{e:c3.2}. Therefore  we have the following theorem.

 \bt\label{t:3.5}Suppose  $s$ is semisimple.
 Then there exist $B_1,B_2>0$
such that
$$
\left\|\frac{x_i(t)}{\|X(t)\|}-\xi_i\right\|\leq B_1e^{-\de t},\hs t\geq1,
$$
and
$$\left|\frac{1}{\|X(t)\|^2}\,{\sum_{i\in J}\langle A x_i(t),\, x_i(t)\rangle}-s\right|\leq B_2e^{-\de t},\hs t\geq 1.
$$
\et
\br\label{r:3.5}
Recall that in the semisimple case we have $\mb{\em GE}_s(A)=\mb{\em E}_s(A)$. Therefore by Lemma \ref{l:3.1} one concludes that
$$
\mb{\em span}\{y_i:\,\,i\in J\}=\mb{\em E}_s(A).
$$
\er

\subsection{A generalized  Perron-like theorem for nonnegative matrices}\label{s:3.3}
\vs
As a corollary of Theorem \ref{t:3.3},  we now state  a generalized  Perron-like theorem for nonnegative matrices.
 Let
$$\R^m_+=\{x=(x_1,\cdots,x_m)\in\R^m:\,\,x_i\geq 0\mb{ for }i\in J\}.$$
A vector $x\in \R^m$ is said to be  {\em nonnegative},
this means that  $x\in \R^m_+$.

\bt\label{t:3.6}Let  $A\in\bbM_m$ be a nonnegative  matrix. Then $r:=r(A)\in\sig(A)$, and  the  principal dominant eigenspace  \,$\mb{\em DE}_r (A)$ has a basis consisting of nonnegative eigenvectors.

In particular, if $r$ is semisimple, then
the principal eigenspace $\mb{\em E}_r(A)$ has a basis consisting of nonnegative eigenvectors.
\et

To prove Theorem \ref{t:3.6}, we need  a well known  positivity property of exponential functions of nonnegative matrices.

\bl\label{l:2.5}Let $A\in \bbM_m$ be nonnegative. Then for any $a\in \R$,
\be\label{e:2.12}
e^{t(A+aI)}\R^m_+\subset \R^m_+,\hs t\geq 0.
\ee
\el
{\bf Proof.} We include a simple proof for the reader's convenience.

Let $P=\R^m_+$. Since $e^{t(A+aI)}P=e^{at}\(e^{tA}P\)$ and $e^{at}>0$, to prove \eqref{e:2.12}, it suffices to verify that $e^{tA}P\subset P$.
 So let $x\in P$, and $t\geq0$. By the nonnegativity of $A$ we see that $A^kx\in P$ (and hence $\frac{t^k}{k!}A^kx\in P$) for all $k\in \bbN$.  Consequently
$
\sum_{k=0}^n\frac{t^k}{k!}A^kx\in P
$
for $n\in\bbN$. Since $P$ is closed, setting $n\ra\8$  one immediately concludes that $e^{tA}x\in P$. $\bx$

\Vs
\noindent{\bf Proof of Theorem \ref{t:3.6}}.
We infer from  Example 2.2 that $A$ is a Perron-like matrix with  $s:=s(A)=r(A):=r$. To prove the theorem, there remains to check that $\mb{DE}_s (A)$ has a basis consisting of nonnegative eigenvectors.

For this purpose, we take
$$
v_i=(\,\stac{i-1}{\overbrace{0,\cdots,0}},1,0,\cdots,0)^{\mb{\tiny T}},\hs i\in J.
$$
 By Lemma \ref{l:2.5} we have $x_i(t)=e^{tA}v_i\in P:=\R^m_+$ for all $t\geq 0$. Consequently
$\frac{x_i(t)}{\|X(t)\|}\in P$ for $t\geq 0$. Let $\xi_i\,\,(i\in J)$ be the vectors given by Theorem \ref{t:3.3}. Then by \eqref{e:c3.1} we deduce that
$\{\xi_i:\,\,\,i\in J\}\subset P.$ This and \eqref{e:c3.20} complete the proof of what we desired. $\bx$

\br As a corollary of Theorem \ref{t:3.6}, we deduce that the principal eigenspace $\mb{\em E}_s(A)$ of a  nonnegative symmetric matrix $A$ has a basis consisting of nonnegative eigenvectors.
\er

\section{Computing  Principal Eigenvalues and Eigenvectors  via Polynomial Approximations of Matrix Exponentials}\label{s:4}
\vs
Based on the convergence results in Section \ref{s:3}, we now show that the principal eigenvalues and eigenvectors of Perron-like matrices   can be  computed via polynomial approximations of matrix exponentials.

We employ  the same notations as in Section \ref{s:3}.

Let $A\in \bbM_m$ be a Perron-like matrix with the principal eigenvalue  $s:=s(A)$, and   $\{v_i\}_{i\in J}$ a basis of $E:=\R^m$. Let $y_i$, $x_i(t)$, $V$, $Y$ and $X(t)$ be the same as in Section \ref{s:3}, and $\xi_i$ ($i\in J$)  the vectors given in Theorem \ref{t:3.3}.
Set
 \be\label{e:Q}\Xi=[\xi_i]_{i\in J},\hs q(t)=[q_i(t)]_{i\in J},\ee  where $q_i(t)=\frac{Q_{y_i}(t)}{\left\|Q_Y(t)\right\|}.$
Then the convergence results in \eqref{e3.1}, \eqref{e:c3.1} and \eqref{e:c3.2} can be reformulated  as
\be\label{e3.0m}
\left\|\frac{X(t)}{\|X(t)\|}-q(t)\right\|\leq B_0 e^{-\de t},\hs t\geq1,
\ee
\be\label{e3.1m}
\left\|\frac{X(t)}{\|X(t)\|}-\Xi\right\|\leq O(t^{-1})+B_1e^{-\de t},\hs t\geq1,
\ee
and
\be\label{e3.2m}\left|\frac{1}{\|X(t)\|^2}\langle AX(t),X(t)\rangle-s\right|\leq O(t^{-1})+ B_2e^{-\de t},\hs t\geq 1.
\ee

\vs
Consider the polynomial approximations of $X(t)$:
$$
X_n(t)=\sum_{k=0}^n\frac{t^k}{k!}A^kV,\hs n\in\bbN.
$$
For computational convenience, in what follows we put $\|V\|=1$.

 Let $\kappa$ be the positive constant in \eqref{e3.28}, and write  \be\label{e:vt}
 \kappa^{-1}e^{2\|A\| t}\(\frac{e\|A\|t}{n+1}\)^{n+1}:=\vartheta_n(t).
 \ee

\bt\label{t:4.0} The following estimates hold true for all $n\in\bbN$:
\be\label{e:4.0}
\left\|\frac{X_n(t)}{\|X_n(t)\|}-q(t)\right\|\leq\,\,B_0e^{-\de t}+\vartheta_n(t)\,o(1),\hs t\geq 1,\ee
\be\label{e4.3}
\left\|\frac{X_n(t)}{\|X_n(t)\| }-\Xi\right\|\leq\,\, O(t^{-1})+B_1e^{-\de t}+\vartheta_n(t)\,o(1),\hs t\geq 1,\ee
and
\be\label{e4.4}\left|\frac{1}{\|X_n(t)\|^2}\langle AX_n(t),X_n(t)\rangle-s\,\right|
\leq  O(t^{-1})+B_2e^{-\de t}+\|A\|\vartheta_n(t)\,o(1),\hs t\geq 1,\ee
where
  $o\(1\)$ is  an infinitesimal  as $n\ra\8$ which  depends only upon $n$.\et
{\bf Proof.}  We observe that
$$
X(t)=\sum_{k=0}^\8\frac{t^k}{k!}A^kV=X_n(t)+R_n(t)V,
$$
where $R_n(t)=\sum_{k=n+1}^\8\frac{t^k}{k!}A^k.$\, Note that
\be\label{e:3.19}\ba{ll}
\left\| R_n(t)V\right\| &\leq \sum_{k=n+1}^\8\frac{t^k}{k!}\left\| A^k V\right\| \leq \sum_{k=n+1}^\8\frac{t^k}{k!}\|A\|^k \\[2ex]
&\leq \frac{(t\|A\|)^{n+1}}{(n+1)!}\,e^{t\|A\|},\hs t\geq 0.\ea
\ee
If we write $\left\| X(t)\right\| =\left\| X_n(t)\right\| +r_n(t)$, then
\be\label{e:3.4}
|r_n(t)|\leq \left\| X(t)-X_n(t)\right\| =\left\| R_n(t)V\right\|.\ee
Let
$$\~R_n(t)V=e^{-s t}t^{-(\nu-1)}R_n(t)V,\hs \~r_n(t)=e^{-s t}t^{-(\nu-1)}r_n(t),$$ where $\nu:=\mb{Cord}\,\(\mb{GE}_s(A)\)$ is the cyclic order  of $\mb{GE}_s(A)$. Then
\be\label{e:5.12}\ba{ll}
|\~r_n(t)|\leq \left\|\~R_n(t)V\right\| &=e^{-s t}t^{-(\nu-1)} \left\| R_n(t)V\right\| \\[2ex]
&\leq \frac{(t\|A\|)^{n+1}}{(n+1)!}\,e^{(\|A\|-s) t}\\[2ex]
&\leq \frac{(t\|A\|)^{n+1}}{(n+1)!}\,e^{2\|A\| t},\hs t\geq 1.
\ea
\ee
Here we have used the simple fact that
$$|s|\leq r(A)\leq \|A\|.$$

Set $$\~X(t)=e^{-st}t^{-(\nu-1)}X(t),\hs \~X_n(t)=e^{-st}t^{-(\nu-1)} X_n(t).$$
Then
\be\label{e:3.5}\ba{ll}
&\,\,\,\,\,\, \left\| \frac{X(t)}{\|X(t)\| }-\frac{X_n(t)}{\left\| X_n(t)\right\| }\right\| =\left\| \frac{\~X(t)}{\left\| \~X(t)\right\| }-\frac{\~X_n(t)}{\left\|\~X_n(t)\right\| }\right\|  \\[3ex]
&=\frac{1}{\left\|\~X(t)\right\| \,\left\| \~X_n(t)\right\| }\left\| {\|\~X_n(t)\| \~X(t)-\|\~X(t)\| \~X_n(t)}\right\|  \\[3ex]
&=\frac{1}{\left\|\~X(t)\right\| \,\left\|\~X_n(t)\right\| }\left\| {\|\~X_n(t)\| \~R_n(t)V-\~r_n(t)\~X_n(t)}\right\|  \\[3ex]
&\leq(\mb{by }\eqref{e:3.4})\leq  \frac{2\left\|\~R_n(t)V\right\| }{\left\|\~X(t)\right\| }.
\ea
\ee
By \eqref{e3.28} we have
$
\|\~X(t)\| \geq \kappa>0$ for $t\geq 1.$
Thus by  \eqref{e:5.12} and \eqref{e:3.5} one concludes that
\be\label{e:3.2}
\left\|\frac{X(t)}{\|X(t)\| }-\frac{X_n(t)}{\|X_n(t)\| }\right\| \leq \frac{2}{\kappa}\frac{(t\|A\|)^{n+1}}{(n+1)!}\,e^{2\|A\| t},\hs t\geq 1.
\ee
Combining this with \eqref{e3.0m} we arrive at the following estimates:
\be\label{e:3.17'}
\left\|\frac{X_n(t)}{\|X_n(t)\| }-q(t)\right\| \leq B_0 e^{-\de t}+\frac{2}{\kappa}\frac{(t\|A\|)^{n+1}}{(n+1)!}\,e^{2\|A\| t},\hs t\geq 1,\ee
and
\be\label{e:3.17}
\left\|\frac{X_n(t)}{\|X_n(t)\| }-\Xi\right\| \leq O(t^{-1})+B_1 e^{-\de t}+\frac{2}{\kappa}\frac{(t\|A\|)^{n+1}}{(n+1)!}\,e^{2\|A\| t},\hs t\geq 1.\ee
Thanks to  the classical Stirling's formula, we have
$$
\lim_{n\ra\8}\frac{n!}{n^{(n+1/2)}e^{-n}}=\sqrt{2\pi},
$$
by which one finds  that
$$
\frac{(t\|A\|)^{n+1}}{(n+1)!}=\(\frac{e\|A\|t}{n+1}\)^{n+1}o(1),
$$
where $o(1)$ is an infinitesimal as $n\ra\8$ depending only upon $n$. The first and second  estimates in the theorem now follow from \eqref{e:3.17'} and  \eqref{e:3.17}.

Using \eqref{e4.3} and a similar argument as in the verification  of \eqref{e:c3.2}, one can obtain  the third  estimate in Theorem \ref{t:4.0}. We omit the details. $\bx$
\Vs

In the case where $s$ is semisimple, by virtue of Theorem \ref{t:3.5} one  can actually remove  the terms $O(t^{-1})$ in the righthand sides of \eqref{e4.3} and \eqref{e4.4}. In other words,  we have a strengthened version of Theorem \ref{t:4.0}:

\bt\label{t:4.4} Assume that $s$ is semisimple. Then
$$
\left\|\frac{X_n(t)}{\|X_n(t)\| }-\frac{Y}{\|Y\|}\right\| \leq B_1e^{-\de t}+\,\vartheta_n(t)o(1),\hs t\geq 1,
$$
 and
$$
\left|\frac{1}{\|X_n(t)\|^2}\left\langle AX_n(t),X_n(t)\right\rangle-s\,\right|\leq B_2e^{-\de t}+\|A\|\vartheta_n(t)o(1),\hs t\geq 1.
$$
\et

\section{An Iterative  Scheme}\label{s:5}
\vs
Theoretically the convergence  results in Sections \ref{s:3} and \ref{s:4} already provide a dynamical way    for computing the principal eigenvalue  and the principal dominant  eigenvectors  of a matrix. However, there is the danger of overflow in calculating $X(t)$ and its polynomial approximations $X_n(t)$ when $t$ and $n$ are  large enough. To overcome  this drawback,  we develop an iterative scheme with double indices corresponding to $t$ and $n$ which  has a more rapid exponential convergence rate than the one given  in Theorem \ref{t:4.0}.

\subsection{An  iterative  scheme}\label{s:5.1}
\vs
 We continue the argument in Section \ref{s:4}. Let $A\in \bbM_m$ be a Perron-like  matrix, and set $T=e^A$. Consider the polynomial approximations of $T$: $$T_n=\sum_{k=0}^n\frac{1}{k!}A^k,\hs n\in \bbN.
$$
For each $n\in \bbN$, we have
$$T=T_n+R_n,\hs\mb{where }\,R_n=\sum_{k=n+1}^\8\frac{1}{k!}A^k.$$
Similar calculations as in \eqref{e:3.19} yield
$$
\|R_n M\| \leq \(\frac{\|A\|^{n+1}}{(n+1)!}\,e^{\|A\|}\)\|M\|,\hs \A\,M\in \bbM_m.
$$

\vs
Let  $$\Sigma_1=\{M\in\bbM_m:\,\,\|M\| =1\}.$$ Given $M\in\Sigma_1$,  we write
$$\|TM\| =\|T_n M\| +r_n(M).$$ Then
\be\label{e:4.2}\ba{ll}
|r_n(M)|&\leq \left|\|TM\|-\|T_n M\| \right|\leq \|TM-T_nM\|\\[2ex]
&= \|R_n M\| \leq \frac{\|A\|^{n+1}}{(n+1)!}\,e^{\|A\|}.
\ea
\ee

Define two mappings $K$ and $K_n$ on $\Sigma_1$ as follows: \,$\A\,M\in \Sigma_1$,
$$K M=\frac{TM}{\|TM\| },\hs K_n M=\frac{T_n M}{\|T_n M\| }.$$  Using some similar calculations as in \eqref{e:3.5} one can  obtain  that
\be\label{e:5.a2}
\|KM-K_nM\| =\left\|\frac{TM}{\|TM\| }-\frac{T_nM}{\|T_nM\| } \right\|\leq \frac{2\|R_nM\| }{\|TM\| },\hs M\in\Sigma_1.
\ee
Observing  that
$$
\|M\| =\|T^{-1}TM\| \leq \|T^{-1}\|\|TM\|,\hs M\in\Sigma_1,
$$
we deduce that
\be\label{e:4.1}
\|TM\| \geq \|T^{-1}\|^{-1}\|M\| \geq e^{-\|A\|},\hs \,M\in\Sigma_1.
\ee
Thus by  \eqref{e:4.2} and \eqref{e:5.a2} it follows that
\be\label{e:5.a5}
\|KM-K_nM\|\leq \lam\,\frac{\|A\|^{n+1}}{(n+1)!}:=o_n,\hs n\in\bbN, 
\ee
where \be\label{e:lam}\lam=2e^{2\|A\|}\geq 2.\ee

\vs Now let   $V=[v_i]_{i\in J}\in \Sigma_1$ be a  nonsingular matrix.
Define an iteration sequence  with double indices as below:
\be\label{e5.6}
M_n(0)=V,\,\,\,\,M_n(k+1)=K_nM_n(k),\hs\,\, k,n\in\bbN.
\ee
Let $q(t)$ and $\Xi$ be the same  as in Theorem \ref{t:4.0}.
\bt\label{t:5.1} Let \,$W_n=M_n(n)$. Then  for $n\in\bbN$, we have
\be\label{e:5.6}
\ba{ll}
\|W_n-q(n)\| \leq B_0e^{-\de n}+\(\frac{\lam\|A\|}{n+1}\)^{n+1}o(1),
\ea
\ee
\be\label{e:5.7}
\ba{ll}
\|W_n-\Xi\| \leq O(n^{-1})+B_1e^{-\de n}+\(\frac{\lam\|A\|}{n+1}\)^{n+1}o(1),
\ea
\ee
 and
 \be\label{e:5.8}\ba{ll}\left|\langle A W_n,\,W_n\rangle -s\,\right|\leq O(n^{-1})+B_2e^{-\de n}+ {\|A\|}\,\(\frac{\lam\|A\|}{n+1}\)^{n+1}o(1),
\ea\ee
 where $o(1)$ is an infinitesimal as $n\ra\8$ depending only upon $n$.
\et
{\bf Proof.} Define a sequence $\{M(k)\}_{k\in\bbN}$ as follows:
$$M(0)=V,\hs M(k+1)=KM(k)\,\,\,(k\in\bbN).$$
Then
$$\ba{ll}
M(k)&=\frac{TM(k-1)}{\|TM(k-1)\| }
=\frac{T\(\frac{TM(k-2)}{\|TM(k-2)\| }\)}{\left\|T\(\frac{TM(k-2)}{\|TM(k-2)\| }\)\right\| }\\[2ex]
&=\frac{T^2M(k-2)}{\|T^2M(k-2)\| }=\cdots=\frac{T^k M_0}{\|T^k M_0\| }=\frac{e^{kA}V}{\|e^{kA}V\| }\,.
\ea$$
Hence by virtue of   \eqref{e3.0m} and \eqref{e3.1m}  we have
\be\label{e:5.9}
\left\|M(k)-q(k)\right\| \leq B_0e^{-\de k}
\ee
and \be\label{e5.7}
\left\|M(k)-\Xi\right\| \leq O(k^{-1})+B_1e^{-\de k}.
\ee

In the sequel we give an  estimate for $\|M_n(k)-M(k)\| $. For notational simplicity, we rewrite $M_n(k)=\widetilde{M}(k)$. Then
\be\label{e:5.a6}\ba{ll}
&\,\,\,\,\,\,\|\widetilde{M}(k)-M(k)\| \\[2ex]
&=\|K_n\widetilde{M}(k-1)-KM(k-1)\| \\[2ex]
&\leq \|K_n\widetilde{M}(k-1)-K\widetilde{M}(k-1)\| +\|K\widetilde{M}(k-1)-KM(k-1)\| \\[2ex]
&\leq(\mb{by }\eqref{e:5.a5})\leq o_n+\|K\widetilde{M}(k-1)-KM(k-1)\|,
\ea\ee
where $o_n$ is the same as  in \eqref{e:5.a5}.
On the other hand,
$$\ba{ll}
\|KX-KY\| &= \left\|\frac{\|TY\| (TX-TY)}{\|TX\| \|TY\| }+ \frac{(\|TY\| -\|TX\| )TY}{\|TX\| \|TY\| }\right\| \\[2ex]
&\leq\frac{\|TX-TY\| }{\|TX\| }+ \frac{|\|TY\| -\|TX\| |}{\|TX\| }\\[2ex]
&\leq \frac{2\|TX-TY\| }{\|TX\| }\leq  \frac{2\|T\|}{\|TX\| } \|X-Y\| \\[2ex]
&\leq(\mb{by }\eqref{e:4.1})\leq  \lam  \|X-Y\|,\hs\A\, X,Y\in\Sigma_1.
\ea
$$
 Hence by \eqref{e:5.a6} one has
$$\ba{ll}
\|\widetilde{M}(k)-M(k)\| &\leq o_n+\lam\|\widetilde{M}(k-1)-M(k-1)\| \leq\cdots\\[2ex]
&\leq o_n(1+\lam+\lam^2+\cdots+\lam^{k-2})+\lam^{k-1}\|\widetilde{M}(1)-M(1)\|.\ea
$$
Because  $\widetilde{M}(0)=M_n(0)=V=M(0)$,  by  \eqref{e:5.a5} we deduce that
$$
\|\widetilde{M}(1)-M(1)\| =\|K_n\widetilde{M}(0)-KM(0)\| =\|K_n V-K V\| \leq o_n.
$$
Thereby
\be\label{e:4.9}
\|\widetilde{M}(k)-M(k)\| \leq o_n(1+\lam+\cdots+\lam^{k-1})= \frac{\lam^k-1}{\lam-1}\,o_n\leq \lam^k \,o_n.
\ee

Combining \eqref{e:5.9}, \eqref{e5.7} and \eqref{e:4.9} it yields
$$\ba{ll}
\|M_n(k)-q(k)\|
\leq B_0e^{-\de k}+{\lam^{k+1}}\frac{\|A\|^{n+1}}{(n+1)!}\,,
\ea $$
and
$$\ba{ll}
\|M_n(k)-\Xi\|\leq O(k^{-1})+B_1e^{-\de k}+{\lam^{k+1}}\frac{\|A\|^{n+1}}{(n+1)!}\,.
\ea $$
Taking $k=n$ and using the Stirling's formula,
we immediately arrive at   the estimates in \eqref{e:5.6} and  \eqref{e:5.7}.

As in Theorem \ref{t:3.3}, the estimate in \eqref{e:5.8} is actually a consequence of \eqref{e:5.7}. We omit the details of the proof.   $\bx$
\Vs

As in Theorem \ref{t:4.4}, if the principal eigenvalue  is semisimple, then  using the convergence result in  Theorem \ref{t:4.4} and repeating the above argument with minor modifications, it can be shown that the iterative scheme developed here has in fact an exponential convergence rate.
More precisely, let $Y=[y_i]_{i\in J}$, where $y_i=\Pi_1 v_i$, and $\Pi_1$ is the projection from $E=\R^m$ to $E_1=\mb{GE}_s(A)$. Write $W_n=M_n(n)$. We have

\bt\label{t:4.1'} Assume that $s$ is semisimple. Then for $n\in\bbN$, we have
\be\label{e:5.7'}
\ba{ll}
\left\|W_n-\frac{Y}{\|Y\|}\right\| \leq B_1e^{-\de n}+\(\frac{\lam\|A\|}{n+1}\)^{n+1}o(1),
\ea
\ee
  and \be\label{e:4.7'}\left|\langle A W_n,\,W_n\rangle -s\,\right|\leq B_2e^{-\de n}+ {\|A\|}\(\frac{\lam\|A\|}{n+1}\)^{n+1}o(1).
\ee
\et
\br\label{r:5.6} Since
$\mb{\em span}\{y_i:\,\,i\in J\}=\mb{\em E}_s(A)$ (see Remark \ref{r:3.5}), in the semisimple case, by \eqref{e:5.7'} the iterative scheme given here  provides an alternative effective  way to compute the whole principal eigenspace $\mb{\em E}_s(A)$ with exponential convergence.
\er

\br Nonnegative irreducible matrices and symmetric matrices are typical examples for which the principal eigenvalues are semisimple. For these two types of matrices one can find a large number  of excellent   works in the literature on the computation of eigenvalues and eigenvectors; see  e.g. \cite{At,Bell, Bra,Buns,Hall,Noda,Parl,Prak,Swar,Wen,WS} etc.\er

\subsection{An iterative scheme with parameters}
\vs
 It can be  seen from Theorems \ref{t:5.1} and \ref{t:4.1'} that  the convergence rate of the scheme can  be significantly affected by  the  scale of the  matrix $A$, and large scale of $A$ may cause slow convergence speed in the early stage of the iteration.
Another risk for a large scale matrix is the overflow in the computation.
 To overcome these  deficiencies, one may introduce a  parameter $\gamma>0$ in the scheme.  Specifically, we  choose an appropriate $\gam$ and use the matrix $\gam A$ in place of $A$. The iterative scheme is then reformulated   as below:

  For each $n\in \bbN$, let
 $$K_n(\gam)X=\frac{T_n(\gam)X}{\|T_n(\gam)X\|},\hs X\in \Sigma_1,$$ where
$
T_n(\gam)=\sum_{k=0}^n\frac{1}{k!}(\gam A)^k.
$
Given a {\em nonsingular} matrix  $V\in \Sigma_1$, set
\be\label{e:4.10}
M_n(0)=V,\hs M_n(k+1)=K_n(\gam)M_n(k)\,\,\,(k\in \bbN).
\ee
Correspondingly   Theorem \ref{t:5.1} now takes the form:

\bt\label{t:5.1'} Let $W_n=M_n(n)$. Then for $n\in\bbN$ we have
$$
\ba{ll}
\|W_n-q(n)\|\leq B_0(\gam)e^{-\gam\de n}+\(\frac{\lam(\gam)\gam\|A\|}{n+1}\)^{n+1}o(1),
\ea
$$
$$
\ba{ll}
\|W_n-\Xi\| \leq O(n^{-1})+B_1(\gam)e^{-\gam\de n}+ \(\frac{\lam(\gam)\gam\|A\|}{n+1}\)^{n+1}o(1),\ea
$$
where $\lam(\gam)=2e^{2\gam\|A\|}$ (see \eqref{e:lam}),  and
 $$\ba{ll}\left|\langle A W_n,\,W_n\rangle -s\,\right|&\leq \gam^{-1}O(n^{-1})+\gam^{-1}B_2(\gam)e^{-\gam\de n}+\\[2ex]
 &\hs+ {\|A\|}\(\frac{\lam(\gam)\gam\|A\|}{n+1}\)^{n+1}o(1).\ea
$$

\et
\br As in Theorem \ref{t:4.1'}, in the case where $s$ is semisimple, the terms $O(n^{-1})$ and $\gam^{-1}O(n^{-1})$ in Theorem \ref{t:5.1'} can be removed. We omit the details.\er
\subsection{A remark on the initial matrix}
\br\label{r:5.5}In the definition of the iteration sequence in \eqref{e5.6}, it is required that $\|V\|=1$. This is just for the sake of convenience in statement.
In practice we can  take a non-normalized  matrix $V$ as the initial one in  the iteration. For such a matrix  one easily examines that the matrices   $M_n(k)$ given by  \eqref{e5.6} (except  $M_n(0)$) and the matrices  $q(n)$ and $\Xi$  in Theorem \ref{t:5.1}  are  the same as in the case where, instead of $V$,  one uses the normalization  $V'=V/\|V\|$ of $V$ as the initial matrix. Hence  all the results above  remain valid.
\er
\subsection{Numerical examples}\label{s:5.2}\vs
We now give several numerical examples to illustrate the efficiency of our scheme.
All the numerical simulations in this work  were carried out on the authors'  PC by using  Matlab. Detailed information concerning the PC and the softwares  is   as follows.
\Vs
\noindent{\bf $\bullet$ Computer Configuration}
\vs
 Processor: Intel(R) Core(TM) i5-4200M; CPU @ 2.50GHZ 2.49GHZ
\vs
 RAM: 8GB
 \vs
 Operation system:  Windows 10 Professional Edition
\vs
 System type: 64-bit operating system, x64-based processor
\Vs
\noindent{\bf $\bullet$ Softwares}\vs
Matlab version: Matlab R2017b

\Vs
\noindent{\bf $\bullet$  Notations}\vs
Given a Perron-like matrix $A$ and an initial matrix $V$ ($V$ is nonsingular), let $\Xi$ and $W_n:=M_n(n)$ be given by Theorem \ref{t:5.1}. Then $W_n$ and $s_n:=\langle A W_n,\,W_n\rangle$ can be regarded as  numerical values of $\Xi$ and $s$, respectively. We will employ the following notations:
$$\ba{ll}
&\ba{ll}\epsilon\(W_n,\Xi\): \,\,\,&\mb{Error estimate  between  $W_n$ and $\Xi$}\ea
\\[2ex]
&\,\,\,\ba{ll}\epsilon(s_n,s): &\,\,\mb{Error estimate  between  $s_n$ and $s$}\ea
\\[2ex]
&\hs\hs\,\,\ba{ll}\epsilon: &\,\,\mb{The total error estimate $\ve(W_n,\Xi)+\ve(s_n,s)$}\ea
\ea
$$
\vs

\noindent{\bf Example 5.1.} Let
$$A=\(\begin{matrix}14&-\frac{21}{2}&\frac{39}{2}&\frac{63}{2}&12&-12&\frac{9}{2}\\[1ex]
                   \frac{18}{5}&-\frac{12}{5}&\frac{28}{5}&\frac{46}{5}&\frac{31}{10}&-\frac{21}{10}&\frac{11}{10}\\[1ex]
                   \frac{12}{5}&\frac{2}{5}&\frac{32}{5}&\frac{34}{5}&\frac{17}{5}&-\frac{27}{5}&\frac{7}{5}\\[1ex]
                   -\frac{26}{5}&\frac{33}{10}&-\frac{87}{10}&-\frac{119}{10}&-\frac{57}{10}&\frac{67}{10}&-\frac{11}{5}\\[1ex]
                   -\frac{8}{5}&\frac{7}{5}&-\frac{13}{5}&-\frac{21}{5}&\frac{2}{5}&\frac{8}{5}&-\frac{3}{5}\\[1ex]
                   \frac{16}{5}&-\frac{3}{10}&\frac{57}{10}&\frac{89}{10}&\frac{21}{5}&-\frac{21}{5}&\frac{17}{10}\\[1ex]
                   \frac{36}{5}&-\frac{63}{10}&\frac{117}{10}&\frac{189}{10}&\frac{36}{5}&-\frac{36}{5}&\frac{47}{10}\end{matrix}\).$$
 The principal eigenvalue $s$ of the matrix is $s= 2$, which is a semisimple  eigenvalue with geometric multiplicity $5$.  Theorem \ref{t:5.1'} indicates that the iteration given in \eqref{e5.6} has an exponential convergence rate for any nonsingular initial matrix $V$. We also infer from \eqref{e:c3.20} and \eqref{e:2sm} that $$\mb{E}_s(A)=\mb{span}\{\xi_i:\,\,i\in J:=\{1,2,\cdots,7\}\},$$ where $\xi_i$ ($i\in J$) are the column vectors of $\Xi$.

 Taking the initial matrix $V$ to be the identity matrix $I$ and using the iterative scheme  in \eqref{e5.6}, we obtain the following table concerning the numerical results by  using the iterative method  developed in this section.

\begin{table}[H]
\caption{Numerical results via the iterative scheme  in \eqref{e5.6} }\label{table:5.1}
\begin{center}
\begin{tabular}{|c|>{\centering\arraybackslash}p{2.5cm}|>{\centering\arraybackslash}p{2.5cm}|>{\centering\arraybackslash}p{2.5cm}|}
\hline
$n$& $\epsilon(W_n,\Xi)$&$\epsilon(s_n,s)$&$\epsilon$\\[1ex]
\hline
1&0.0785&0.0980&0.1765\\
\hline
2&0.0076&0.0371&0.0447\\
\hline
3&$7.0861\times10^{-5}$&$3.7401\times10^{-4}$&$4.4487\times10^{-4}$\\
\hline
4&$1.2770\times10^{-6}$&$2.4147\times10^{-5}$&$2.5424\times10^{-5}$\\
\hline
5&$8.6196\times10^{-8}$&$9.1289\times10^{-7}$&$9.9909\times10^{-7}$\\
\hline
8&$9.9724\times10^{-12}$&$1.0554\times10^{-10}$&$1.1552\times10^{-10}$\\
\hline
10&$5.0047\times10^{-14}$&$2.4603\times10^{-13}$&$2.9607\times10^{-13}$\\
\hline
\end{tabular}
\end{center}
\end{table}

\Vs
\noindent{\bf Example 5.2.}
Consider the  symmetric matrix
$$A=\(\begin{matrix}1&1&1&1\\
                   1&1&-1&-1\\
                   1&-1&1&-1\\
                   1&-1&-1&1\end{matrix}\),$$
for which  $s=s(A)=2$ with algebraic and geometric multiplicity $3$. Taking
$$V=\(\begin{matrix}1&1&1&0\\
                   1&0&0&1\\
                   0&1&0&1\\
                   0&0&1&0\end{matrix}\)$$
and applying the iterative scheme  in \eqref{e5.6}, we obtain the following table concerning the error estimates of the numerical results:

\begin{table}[H]
\caption{Numerical results via the iterative scheme  in \eqref{e5.6}}\label{table:5.2}
\begin{center}
\begin{tabular}{|c|>{\centering\arraybackslash}p{2.5cm}|>{\centering\arraybackslash}p{2.5cm}|>{\centering\arraybackslash}p{2.5cm}|}
\hline
$n$& $\epsilon(W_n,\Xi)$&$\epsilon(s_n,s)$&$\varepsilon$\\[1ex]
\hline
1&0.1252&0.0625&0.1877\\
\hline
2&0.0151&$9.1408\times10^{-4}$&0.0160\\
\hline
3&$5.5105\times10^{-5}$&$1.2146\times10^{-8}$&$5.5117\times10^{-5}$\\
\hline
4&$1.9435\times10^{-6}$&$1.5108\times10^{-11}$&$1.9435\times10^{-6}$\\
\hline
5&$2.4565\times10^{-11}$&$2.2204\times10^{-16}$&$2.4565\times10^{-11}$\\
\hline
10&$1.4687\times10^{-16}$&0&$1.4687\times10^{-16}$\\
\hline
\end{tabular}
\end{center}
\end{table}

The principal eigenvalues of the matrices in the above two examples are semisimple. Numerical results also indicate that in such a case the convergence rate of the iterative scheme developed in this section can be exponential. However, situations seem to be quite different in the case where the principal eigenvalue is not semisimple, and both the theoretical and numerical results demonstrate that the convergence can be very slow.
\Vs
\noindent{\bf Example 5.3.} Consider the matrix
$$A=\(\begin{matrix}2&1&0&0&2\\
                   0&2&1&0&0\\
                   0&0&2&0&1\\
                   0&0&0&1&0\\
                   0&0&0&3&1 \end{matrix}\).$$
The principal  eigenvalue is  $s=2$.  The algebraic multiplicity of $s$ is $3$, whereas it has a geometric multiplicity $1$. Taking  $V=I$ and applying  the scheme given in Theorem \ref{t:5.1} to compute the numerical values $s_n$ of $s$, we obtain the following table, from which it can be easily seen that the numerical results and the theoretical one in \eqref{e:5.8} coincide.

\begin{table}[H]
\caption{Numerical results via  the scheme in  Theorem \ref{t:5.1}}\label{table:6.1}
\begin{center}
\begin{tabular}{|c|>{\centering\arraybackslash}p{2.5cm}|>{\centering\arraybackslash}p{2.5cm}|}
\hline
$n$&$s_n$&$\epsilon(s_n,s)$\\[2ex]
\hline
10&2.2160&0.2160\\
\hline
50&2.0413&0.0413\\
\hline
100&2.0203&0.0203\\
\hline
200&2.0101&0.0101\\
\hline
500&2.0040&0.0040\\
\hline
\end{tabular}
\end{center}
\end{table}

\section{Determination of the Cyclic Order  of $\mb{GE}_s(A)$}\label{s:6}
\vs
We have seen from either the theoretical results in Theorems \ref{t:3.3}, \ref{t:4.0} and \ref{t:5.1} or the numerical example in Example 5.3 that the situation is more complicated and worse in the case where the principal eigenvalue $s$ of a matrix $A$ is non-semisimple. On one hand, the computation methods presented in the previous sections can only allow us to obtain numerical results for the principal dominant eigenspace $\mb{DE}_s(A)$ (recall that all the column vectors of the matrix $\Xi$ in the aforementioned theorems are principal dominant eigenvectors). On the other hand, the convergence speed  can be very slow, which  fact may make the methods   to be of little practical  sense.

The remaining part of this paper is devoted to this bad case  where $s$ is non-semisimple. Our main aim is to  present  a combined  method for the computation of the principal eigenvalue and develop an efficient scheme to compute the whole principal generalized eigenspace $\mb{GE}_s(A)$.
To this end, we first  address the problem  of how to determine the cyclic order $\mb{Cord}\(\mb{GE}_s(A)\)$ of the space $\mb{GE}_s(A)$, which is also of  independent interest in its own right.

\subsection{Determining  the cyclic order  of $\mb{\em GE}_s(A)$}\vs
We use the same notations as in Section \ref{s:3}.
So let $v_i$, $y_i$ and $z_i$ ($i\in J$) be the same as therein.
Let $Q_{y_i}(t)$ be the characteristic polynomial of $y_i$, and $Q_Y(t)=\[Q_{y_i}(t)\]_{i\in J}$. We infer from the proof of Theorem \ref{t:3.1} that
\be\label{e6.1}
t^{-(\nu-1)}\left\|Q_Y(t)\right\|=\omega+O(t^{-1}),
\ee
where $\omega$ is the number in \eqref{e3.19}.
Let $q_i(t)$ and $q(t)$ be given as in \eqref{e:Q}.

Assume that $j\in J$ is a number  such  that $y_j\ne 0$.  For each $k\in\{0,1,\cdots,\nu\}$, where $\nu=\mb{Ord}(y_j)$, put
 \be\label{e:6.38}
 \psi_k(t,\bar s)=t^{2k}{\|(A-\bar s I)^k q_j(t)\|^2},
 \ee
where $\bar s:=s+\De s$ denotes  an approximate value of $s$.
Note that
\be\label{e6.3}
(A-\bar s I)^k=\sum_{i=0}^kC_k^i(-\De s)^i(A-s I)^{k-i},
\ee
where $C_k^i$ are the combinatorial numbers. We  observe that
\be\label{e6.4}\ba{ll}
(A-s I)^{k} Q_{y_j}(t)&=\sum_{\ell=0}^{\nu-1}\frac{t^\ell}{\ell!}(A-s I)^{k+\ell} y_j\\[2ex]
&=\sum_{\ell=0}^{\nu-k-1}\frac{t^\ell}{\ell!}(A-s I)^{k+\ell} y_j\\[2ex]
&=t^{\nu-k-1}\(\frac{1}{(\nu-k-1)!}\~w_j+O(t^{-1})\)
\ea\ee
for $k\leq \nu-1$, where  \be\label{e:6.2}\~w_j=(A-sI)^{\nu-1}y_j.\ee

\vs In what follows we assume that  $|t\De s|\leq 1$ and  argue by cases.
\Vs
{\bf Case 1}: \,$k\leq\nu-1$.\, In this case by \eqref{e6.3} and \eqref{e6.4} we deduce   that
 \be\label{e:6.20}\ba{ll}
(A-\bar s I)^k Q_{y_j}(t)&=\sum_{i=0}^kC_k^i(-\De s)^i  t^{\nu-k+i-1}\(\frac{\~w_j}{(\nu-k+i-1)!}+O(t^{-1})\)\\[2ex]
&=t^{\nu-k-1}\sum_{i=0}^k(-t\De s)^i\(a_{k,i}\~w_j+O(t^{-1})\)\\[2ex]
&=t^{\nu-k-1}a_{k,0}\~w_j+ t^{\nu-k-1}\sum_{i=1}^k(-t\De s)^i\,a_{k,i}\~w_j+\\[2ex]
&\hs + t^{\nu-k-1}\sum_{i=0}^k(-t\De s)^iO(t^{-1})\\[2ex]
&=t^{\nu-k-1}\(a_{k,0}\~w_j+ O(|t\De s|)+ O(t^{-1})\)
\ea \ee
as $t\ra\8$ and $t\De s\ra 0$,  where
 \be\label{e:6.35}a_{k,i}=\frac{C_k^i}{(\nu-k+i-1)!}\,.\ee
 Therefore by \eqref{e6.1} one easily  examines that
 \be\label{e6.7}
 t^k(A-\bar s I)^kq_j(t)=\omega^{-1}a_{k,0}\~w_j+O(|t\De s|)+O(t^{-1}).
 \ee
 Hence
  \be\label{e6.8}\ba{ll}
\psi_k(t,\bar s)=\omega^{-2}\|a_{k,0}\~w_j\|^2+O(|t\De s|)+O(t^{-1}).
\ea
 \ee

\vs
{\bf Case 2}: \,$k=\nu. $ Noticing that $(A-s I)^k Q_{y_j}(t)=0$, as in \eqref{e:6.20} we have
 \be\label{e:6.36}\ba{ll}
(A-\bar s I)^k Q_{y_j}(t)&=t^{\nu-k-1}\sum_{i=1}^k(-t\De s)^i\(a_{k,i}\~w_j+O(t^{-1})\)\\[2ex]
&=t^{\nu-k-1} O(|t\De s|),
\ea \ee
 Thus by \eqref{e6.1} we deduce that
 \be\label{e6.10}
 t^k(A-\bar s I)^kq_j(t)=O(|t\De s|).
 \ee
It then follows by the definition of $\psi_\nu$ in  \eqref{e:6.38}  that
  \be\label{e:6.42}\ba{ll}
\psi_\nu(t,\bar s)=O\(|t\De s|^2\).
\ea
 \ee

 In practice the function $q_j(t)$ is in fact  unknown. However, Theorem \ref{t:3.1} indicates that $\bar q_j(t):=\frac{x_j(t)}{\|X(t)\|}$ approaches $q_j(t)$ exponentially, where $x_j(t)$ and $\|X(t)\|$ are the same as in Theorem \ref{t:3.1}. More precisely, we have
 \be\label{e:6.39}
 \|\bar q_j(t)-q_j(t)\|\leq B_0e^{-\de t},\hs t\geq 0.
 \ee
 Set  \be\label{e:6.40}
\bar \psi_k(t,\bar s)=t^{2k}{\|(A-\bar s I)^k \bar q_j(t)\|^2}.
 \ee
Combining \eqref{e6.8}, \eqref{e:6.42} and \eqref{e:6.39} together, we  obtain  the following fundament result:

\bp\label{p:6.1}Assume that $y_j\ne0$. Let $\nu=\mb{\em Ord}(y_j)$.  Then as $t\ra\8$ and $t\De s\ra 0$, the following assertions hold: \be\label{e:6.43}
\bar\psi_k(t,\bar s)=\omega^{-2}\|a_{k,0}\~w_j\|^2+O(|t\De s|)+O(t^{-1})+O(t^{2k}e^{-\de t})
\ee
for $k\leq\nu-1$, where $\~w_j$ is given by \eqref{e:6.2}; and
 \be\label{e:6.44}
\bar\psi_\nu(t,\bar s)=O\(|t\De s|^2\)+O(t^{2\nu}e^{-\de t}).
\ee
\ep

 In applications one can use the $j$-th colum vector $W_{n,j}$ of the matrix
 $$W_n:=M_n(n)$$ given in \eqref{e5.6} to replace $\bar q_j(t)$ in the definition of $\bar \psi_k(t,\bar s)$. That is, we define
\be\label{e:6.45}
\bar\psi_k(n,\bar s)=n^{2k}{\|(A-\bar s I)^k W_{n,j}\|^2},\hs k,n\geq 0.
\ee
Thanks to  Theorem \ref{t:5.1}, Proposition \ref{p:6.1} is  correspondingly modified as below:

\bp\label{p:6.2}Assume that $y_j\ne0$. Let $\nu=\mb{\em Ord}(y_j)$.  Then as $n\ra\8$ and $n\De s\ra 0$, we have
\be\label{e:6.46}\ba{ll}
\bar\psi_k(n,\bar s)&=\omega^{-2}\|a_{k,0}\~w_j\|^2+O(n|\De s|)+O(n^{-1})+O(n^{2k}e^{-\de n})+\\[2ex]
&\hs + n^{2k}\(\frac{\lam\|A\|}{n+1}\)^{n+1}o(1),\hs\,0\leq k\leq\nu-1,
\ea
\ee
 \be\label{e:6.47}
\bar\psi_\nu(n,\bar s)=O\((n|\De s|)^2\)+O(n^{2\nu}e^{-\de n})+n^{2\nu}\(\frac{\lam\|A\|}{n+1}\)^{n+1}o(1),
\ee
where $\lam$ is the  constant given in Theorem \ref{t:5.1}.
\ep

Now we assign $$\bar\psi_{-1}(n,\bar s)=\bar\psi_{0}(n,\bar s),$$ and define
\be\label{e:6.48}
\b_k(n,\bar s)=\frac{\bar\psi_k(n,\bar s)}{\bar\psi_{k-1}(n,\bar s)},\hs 1\leq k\leq \nu.
\ee
Clearly \be\label{e:6.34}\b_0(n,\bar s)=1,\hs n\in\bbN.\ee
If we choose for each $n$ an approximate value $\bar s=s_n$ for  $s$ with $n(s_n-s)\ra 0$ as $n\ra\8$, then by Proposition \ref{p:6.2} it can be easily seen that
\be\label{e:6.50}
\lim_{n\ra \8}\b_\nu(n,s_n)=0,
\ee
whereas \be\label{e:6.49}
\lim_{n\ra \8}\b_k(n,s_n)=\nu-k-1\geq 1,\hs k=1,2,\cdots,\nu-1.
\ee

 These observations form a basis for the computation of the cyclic order of the principal generalized eigenspace.

\Vs
\noindent {\bf $\bullet$ The computation scheme for  $\mb{Cord}\(\mb{GE}_s(A)\)$}.\vs

We are now ready to  propose a scheme to compute  the cyclic order $\mb{Cord}\(\mb{GE}_s(A)\)$ of the principal generalized eigenspace.
\vs
{\bf Step 1}. \,Pick an appropriately  small  number $\ve>0$ (say, $\ve=1/10$).
\vs
{\bf Step 2}. \,Let $W_n:=M_n(n)$, where  $M_n(k)$  is generated   by \eqref{e5.6} with $V=I$. Denote by $W_{n,i}$ the $i$-th column vector of $W_n$.

Take a suitably large  number $N\in\bbN$ and compute $W_N$ and $s_N:=\langle AW_N,\,W_N\rangle$. Theorem \ref{t:5.1} asserts that in general
$|s_N-s|=O(N^{-1})$.

 Fix an index  $j\in J$ with
\be\label{ej}\|W_{N,j}\|=\max_{i\in J}\|W_{N,i}\|.\ee
Since $$\ba{ll}\|W_N\|=\(\sum_{i\in J}\|W_{N,i}\|^2\)^{1/2}=1,\ea$$ we deduce that $\|W_{N,j}\|\geq 1/\sqrt m$. Further by Theorem \ref{t:5.1} we see that the $j$-th column vector $\xi_j$ of\, $\Xi$ is nonzero provided that $N$ is sufficiently large.
Since each nonzero column vector of $\Xi$  is a dominated principal eigenvector (see Theorem \ref{t:3.3}), one concludes that \be\label{ey_j}\mb{Cord}\(\mb{GE}_s(A)\)=\mb{Ord}(y_j).\ee
\vs
{\bf Step 3}. Take a   number $n\in\bbN$ with  $n<<N$ and compute $\b_k(n,s_N)$ for $k\in J$.
  Note  that
\be\label{e:6.51}
\b_k(n,\bar s)=\frac{n^2\|(A-\bar s I)^kW_{n,j}\|^2}{\|(A-\bar s I)^{k-1}W_{n,j}\|^2},\hs k\geq 1\,.
\ee

Suppose that there is a number  $k_0\in J$ such that $\b_k(n,s_N)$ demonstrates  the following dichotomy property:
\be\label{e:6.54}
\b_k(n,s_N)\geq 1-\ve\,\,\,(0\leq k<k_0),\hs \b_{k_0}(n,s_N)<\ve.
\ee
Then we  define $\mb{Cord}\(\mb{GE}_s(A)\)=k_0$, and the procedure  ends.

Otherwise, we re-choose  the numbers  $N$ and $n$ and repeat the above procedure until  the  dichotomy phenomenon  in \eqref{e:6.54} appears.

\br One can fix $j$ and change $N$ and $n$ suitably to see whether the phenomenon in  \eqref{e:6.54} becomes stable. \er

\subsection{A numerical example}\vs
\noindent{\bf Example 6.1.} Let $A$ be the matrix given in Example 5.3.
The principal  eigenvalue is  $s=2$. Since the algebraic multiplicity of $s$ is $3$ whereas its geometric multiplicity is $1$, it is easy to  see that $\mb{Cord}\(\mb{GE}_s(A)\)=3$.
Now suppose that $\mb{Cord}\(\mb{GE}_s(A)\):=\nu$ is unknown. Let us try to  use the method proposed above to compute this algebraic quantity.

\Vs
(1) \,Set $\ve=0.10$.
\vs
Let $V=I$. Taking $N=100,\,105,\,110$ and performing the iteration in \eqref{e5.6}, we obtain respectively that  $s_N=2.0203,\,2.0194,\,2.0185$. The number $j=4$ satisfies  \eqref{ej}.

Computing $\b_k(n,s_N)$ for $k\in J$ with different choices of  $n$ ($n<<N$), we get the following tables:

\begin{table}[H]
\caption{$N=100,\,\,s_N=2.0203,\,\,j=4$}\label{table:6.3}
\begin{center}
\begin{tabular}{|c|>{\centering\arraybackslash}p{2.2cm}|>{\centering\arraybackslash}p{1.5cm}|>{\centering\arraybackslash}p{2cm}|>{\centering\arraybackslash}p{2cm}|>{\centering\arraybackslash}p{1.5cm}|}
\hline
$n$&$\b_1(n,s_N)$&$\b_2(n,s_N)$&$\b_3(n,s_N)$&$\b_4(n,s_N)$&$\b_5(n,s_N)$\\[2ex]
\hline
5&4.0341&1.7642&0.1045&47.2091&11.1716\\
\hline
6&4.1960&1.5385&0.0372&57.1796&22.6585\\
\hline
7&4.2764&1.3296&0.1202&6.2138&41.2745\\
\hline
8&4.2902&1.1541&0.2415&0.7801&83.9131\\
\hline
9&4.2511&1.0068&0.3771&0.0721&249.3744\\
\hline
10&4.1760&0.8821&0.5273&0.0592&79.0894\\
\hline
\end{tabular}
\end{center}
\end{table}

\begin{table}[H]
\caption{$N=105,\,\,s_N=2.0194,\,\,j=4$}
\begin{center}
\begin{tabular}{|c|>{\centering\arraybackslash}p{2.2cm}|>{\centering\arraybackslash}p{1.5cm}|>{\centering\arraybackslash}p{2cm}|>{\centering\arraybackslash}p{2cm}|>{\centering\arraybackslash}p{1.5cm}|}
\hline
$n$&$\b_1(n,s_N)$&$\b_2(n,s_N)$&$\b_3(n,s_N)$&$\b_4(n,s_N)$&$\b_5(n,s_N)$\\[2ex]
\hline
4&3.7754&1.8266&0.8536&10.9190&3.6452\\
\hline
5&4.0519&1.7792&0.1081&44.8683&11.0745\\
\hline
6&4.2179&1.5564&0.0326&64.2748&22.3722\\
\hline
7&4.3022&1.3499&0.1055&7.0344&40.2877\\
\hline
8&4.3199&1.1765&0.2150&0.8945&79.4518\\
\hline
9&4.2843&1.0311&0.3369&0.0849&229.0532\\
\hline
10&4.2127&0.9079&0.4705&0.0490&103.6242\\
\hline
\end{tabular}
\end{center}
\end{table}

\begin{table}[H]
\caption{$N=110,\,\,s_N=2.0185,\,\,j=4$}
\begin{center}
\begin{tabular}{|c|>{\centering\arraybackslash}p{2.2cm}|>{\centering\arraybackslash}p{1.5cm}|>{\centering\arraybackslash}p{2cm}|>{\centering\arraybackslash}p{2cm}|>{\centering\arraybackslash}p{1.5cm}|}
\hline
$n$&$\b_1(n,s_N)$&$\b_2(n,s_N)$&$\b_3(n,s_N)$&$\b_4(n,s_N)$&$\b_5(n,s_N)$\\[2ex]
\hline
4&3.7891&1.8380&0.8602&10.6370&3.6143\\
\hline
5&4.0698&1.7942&0.1119&42.5756&10.9797\\
\hline
6&4.2398&1.5744&0.0286&72.0950&22.0971\\
\hline
7&4.3281&1.3704&0.0919&8.0137&39.3580\\
\hline
8&4.3496&1.1991&0.1904&1.0288&75.3986\\
\hline
9&4.3177&1.0555&0.2996&0.1016& 207.5873\\
\hline
10&4.2494&0.9340&0.4181&0.0404&136.9586\\
\hline
\end{tabular}
\end{center}
\end{table}

We see  that if $n=6$, then
$$
\b_k(n,s_N)>1-\ve\,\,\,(0\leq k\leq 2),\hs \b_3(n,s_N)<\ve
$$
for all $N=100,\,105,\,110$. Hence one concludes that $\nu=3$.
\Vs
(2) \,We can also try with $N=1000$. In such a case we obtain that $s_N=2.0020$ and  $j=4$. Taking different $n<<N$ and computing  $\b_k(n,s_N)$ for $k\in J$, it yields the following table:

\begin{table}[H]
\caption{$N=1000,\,\,s_N=2.0020,\,\,j=4$}\label{table:6.2}
\begin{center}
\begin{tabular}{|c|>{\centering\arraybackslash}p{2.2cm}|>{\centering\arraybackslash}p{1.5cm}|>{\centering\arraybackslash}p{2cm}|>{\centering\arraybackslash}p{2cm}|>{\centering\arraybackslash}p{1.5cm}|}
\hline
$n$&$\b_1(n,s_N)$&$\b_2(n,s_N)$&$\b_3(n,s_N)$&$\b_4(n,s_N)$&$\b_5(n,s_N)$\\[2ex]
\hline
5&4.4054&2.0797&0.2263&14.7423&9.6749\\
\hline
10&4.9526&1.4635&0.0024&6.0700&78.3099\\
\hline
15&4.7560&1.2232&0.0085&0.0012&186.1411\\
\hline
20&4.5468&1.1016&0.0155&0.0065& 0.0029\\
\hline
\end{tabular}
\end{center}
\end{table}
\noindent Clearly for $n=10$, we have
$$
\b_k(n,s_N)>1-\ve\,\,\,(k\leq 2),\hs \b_3(n,s_N)<\ve
$$
with $\ve=0.01$. Hence once again we conclude that $\nu=3$.

\section{A Combined Method for the Computation of  Non-semisimple  Principal  Eigenvalues}\label{s:7}
\vs
  Theorem \ref{t:5.1} asserts that $s_n:=\langle AM_n(n),\,M_n(n)\rangle$ converges to the principal eigenvalue $s$ of $A$. If $s$ is semisimple, we know that the convergence is exponentially fast, and therefore  the numerical values of $s$ within tolerance error can be expected to be  computed quickly. However, in the non-semisimple  case  the convergence  may be very slow,  as is also demonstrated in Example 5.3.

On the other hand,  we have noticed in \eqref{e:5.6} that $W_n:=M_n(n)$ approaches $q(n)$ swiftly. Using this simple fact, in this section we propose a new combined method for the computation of non-semisimple  principal eigenvalues, which turns out to be more efficient.

Since we are only interested in the case of non-semisimple  case in this section,  in what follows we always assume that $$\nu:=\mb{Cord}\(\mb{GE}_s(A)\)\geq 2.$$

We will  employ the same notations in  the preceding sections.
 Let $E:=\R^m$. Denote by $\mB_x(r)$ the ball in $E$ centered at $x\in E$ with radius $r>0$.

\subsection{A fundamental lemma}\vs

We begin our work in this section with a fundamental lemma.
\bl\label{l:7.1}Let $\xi\in \mb{\em GE}_s(A)$, $\xi\ne0$. Write $\nu=\mb{\em Ord}(\xi)$. Given $x\in E$, define a function $\phi_x(\tau)$ on $\R$ as
\be\label{e:7.0}
\phi_x(\tau)=\|(A-\tau I)^{\nu}x\|^2,\hs \tau\in\R.
\ee
Then  there is an $\ve>0$ such that for every  $x\in \ol{\mb{\em B}}_\xi(\ve)$, the function  $\phi_x(\tau)$ has a unique minimum point $\tau_x$. Furthermore, there exist $L_0,L_1>0$ such that
\be\label{e:7.5}
|\tau_x-s|\leq L_0\|x-\xi\|,\hs \A\,x\in \ol{\mb{\em B}}_\xi(\ve)
\ee
and
\be\label{e:7.5b}
\phi_x(\tau_x)\leq L_1\|x-\xi\|,\hs \A\,x\in \ol{\mb{\em B}}_\xi(\ve).
\ee
\el
{\bf Proof.} Since $(A-\tau I)^{\nu}\xi=0$ {if and only if} $\tau=s$, the function  $\phi_\xi(\tau)$ has a unique minimum point $\tau_\xi=s$ with $\phi_\xi(s)=0$. It is a basic knowledge that $\phi_\xi'(s)=0$.
Simple calculations yield
$$\phi_\xi''(\tau)=\left\{\ba{ll}2\|\xi\|^2,\hs &\nu=1\\[1ex]
2\nu(\nu-1) \left\langle(A-\tau I)^{\nu-2}\xi,\,\,(A-\tau I)^{\nu}\xi\right\rangle+
2\nu^2\|(A-\tau I)^{\nu-1}\xi\|^2,\hs&\nu>1. \ea\right.$$
In particular, since $(A-s I)^{\nu-1}\xi\ne 0$, we have
\be\label{e:7.1}\phi_\xi''(s)=2\nu^2\|(A-sI)^{\nu-1}\xi\|^2:=3\b>0.\ee

  Fix an $\eta>0$ such that
$$
\phi_\xi''(\tau)\geq 2\b,\hs |\tau-s|\leq \eta.
$$
Then there is an $\ve_0>0$ such that
\be\label{e:7.15}
\phi_x''(\tau)\geq \b,\hs |\tau-s|\leq \eta,\,\,x\in\ol\mB_\xi(\ve_0).
\ee
We may assume that $\ve_0\leq \|\xi\|/2$. Then it can be easily seen that $\phi_x(\tau)\ra\8$ as $|\tau|\ra\8$ uniformly with respect to $x\in \ol\mB_\xi(\ve_0)$.  Thus one can pick  a $T>0$ such that
\be\label{e:7.16}
\phi_x(\tau)\geq 2,\hs \A\,x\in \ol\mB_\xi(\ve_0),\,\,|\tau|>T.
\ee
It can be assumed that $|s|<T/2$.

Because $\phi_\xi(\tau)>0$ for $\eta\leq |\tau-s|\leq T$, 
we deduce that
$$
\min_{\eta\leq |\tau-s|\leq T}\phi_\xi(\tau):=3c>0.
$$
By continuity  there is $\ve_1>0$ such that
\be\label{e:7.17}
\min_{\eta\leq |\tau-s|\leq T}\phi_x(\tau)\geq 2c,\hs x\in\ol\mB_\xi(\ve_1).
\ee
On the other hand, since  $\phi_\xi(s)=0$,  one can take a positive number $\ve<\min(\ve_0,\ve_1)$ such that \be\label{e:7.18}
\min_{|\tau-s|\leq \eta}\phi_x(\tau)\leq \min(1,c),\hs\A\,x\in \ol\mB_\xi(\ve).
\ee
Combining  \eqref{e:7.16}, \eqref{e:7.17} and \eqref{e:7.18} we deduce that for each $x\in\ol\mB_\xi(\ve)$,
\be\label{e:7.19}
\min_{\tau\in \R}\phi_x(\tau)=\min_{|\tau-s|\leq \eta}\phi_x(\tau).
\ee

\eqref{e:7.15} and \eqref{e:7.19} imply  that for each $x\in\ol\mB_\xi(\ve)$, the function $\phi_x(\tau)$ has a unique minimum point $\tau_x$ with $|\tau_x-s|\leq \eta.$ As $\phi_x'(\tau_x)=0=\phi_\xi'(s)$, we have
$$
\phi_x'(\tau_x)-\phi_x'(s)=\phi_\xi'(s)-\phi_x'(s).
$$
By virtue of  the mean value theorem we deduce that $\phi_x'(\tau_x)-\phi_x'(s)=\phi_x''(\theta)(\tau_x-s)$ for some $\theta$ with $|\theta-s|\leq \eta.$ Therefore by \eqref{e:7.15} we deduce that
\be\label{e:6.21}
\b |\tau_x-s|\leq |\phi_\xi'(s)-\phi_x'(s)|=|\phi_x'(s)|.
\ee
Let  $z=x-\xi$.  We observe that
$$
\ba{ll}
|\phi_x'(s)|&=2\nu\left|\llg(A-s I)^{\nu-1}x,\,(A-s I)^{\nu}x\rrg\right|\\[2ex]
&=2\nu\left|\llg(A-s I)^{\nu-1}x,\,(A-s I)^{\nu}(\xi+z)\rrg\right|\\[2ex]
&=2\nu \left|\llg(A-s I)^{\nu-1}x,\,(A-s I)^{\nu}z\rrg\right|\leq c_0\|z\|.
\ea
$$
Therefore by \eqref{e:6.21} one concludes  that
$
 |\tau_x-s|\leq c_0\b^{-1}\|z\|,
$
which is what we desired in  \eqref{e:7.5}.

We also note that
\be\label{e:7.23}\ba{ll}
\phi_x(\tau_x)
&=\phi_x(s)+\(\phi_x(\tau_x)-\phi_x(s)\)\\[2ex]
&\leq \phi_x(s)+|\phi_x(\tau_x)-\phi_x(s)|.
\ea
\ee
On the other hand, it is easy to verify that
$$
|\phi_x(\tau_x)-\phi_x(s)|\leq L|\tau_x-s|\leq (\mb{by }\eqref{e:7.5})\leq c_1\|z\|
$$
for $x\in\ol\mB_\xi(\ve)$, and
$$\ba{ll}
\phi_x(s)&=\llg(A-s I)^\nu x,\,(A-s I)^\nu x\rrg\\[2ex]
&=\llg(A-s I)^\nu z,\,(A-s I)^\nu z\rrg\leq c_2\|z\|^2.
\ea
$$
Hence by \eqref{e:7.23} we deduce  that
$$
\phi_x(\tau_x)\leq c_1\|z\|+c_2\|z\|^2\leq (c_1+c_2\ve)\|z\|.
$$
This  completes the proof of \eqref{e:7.5b}. $\bx$

\subsection{A combined  method for computing  the principal eigenvalue}

\vs
Now we propose a combined  method for the computation of the principal eigenvalue $s$ of a Perron-like matrix $A$. The basic idea is to pick a suitably large $n$ and solve the minimum point $\tau_*:=\tau_x$ of a suitable   function $\phi_x(\tau)$ as defined by \eqref{e:7.0} with $x=W_{n,j}$,
where $W_{n,j}$ is a  column vector of the matrix $$W_n:=M_n(n)$$ who has the maximal norm among the column vectors $W_{n,i}$ ($i\in J$) of $W_n$ (hence $\|W_{n,j}\|\geq 1/\sqrt m$). For notational convenience, we rewrite $$\phi_{W_{n,j}}(\tau)=\phi_n(\tau).$$

Suppose that $n$ is chosen suitably large, and let  $j\in J$ be  such that $\|W_{n,j}\|=\max_{i\in J}\|W_{n,i}\|$. Then we infer from \eqref{ey_j} that $$\mb{Cord}\(\mb{GE}_s(A)\)=\mb{Ord}(y_j)=\mb{Ord}\(Q_{y_j}(n)\)=\mb{Ord}(\xi):=\nu,$$
where
$
\xi:=q_j(n)=\frac{Q_{y_j}(n)}{\|Q_Y(n)\|}.
$
In view of Theorem \ref{t:5.1} we have
\be\label{e:7.20}
\|W_{n,j}-\xi\| \leq B_0e^{-\de n}+\(\frac{\lam\|A\|}{n+1}\)^{n+1}o(1),
\ee
which indicates that the function $\phi_n(\tau)$ should take the form:
\be\label{e:7.26}\phi_n(\tau)=\|(A-\tau I)^{\nu}W_{n,j}\|^2,\hs \tau\in\R,\ee

Before describing  the combined computation method, let us make a simple observation on the second order derivative  $\phi_n''(\tau)$ near $s$.
 We have
$$\ba{ll}\phi_n''(\tau)&=2\nu(\nu-1) \left\langle(A-\tau I)^{\nu-2}W_{n,j},\,(A-\tau I)^{\nu}W_{n,j}\right\rangle\\[2ex]
&\hs+
2\nu^2\|(A-\tau I)^{\nu-1}W_{n,j}\|^2.\ea$$
Let $\De s=\tau-s$. We also assume that $n|\De s|\leq 1$. Using some similar calculations as in the verifications of \eqref{e6.7}, \eqref{e6.8}, \eqref{e:6.46} and \eqref{e:6.47}, we can get that
$$\ba{ll}
(A-\tau I)^kW_{n,j}&=n^{-k}\(\omega^{-1}a_{k,0}\~w_j+O(n|\De s|)+O(n^{-1})\)+\\[2ex]
&\hs +O(e^{-\de n})+ \(\frac{\lam\|A\|}{n+1}\)^{n+1}o(1)=O(n^{-k})
\ea
$$
for $0\leq k\leq\nu-1$, where $a_{k,0}$ and $\~w_j$ are the same as in \eqref{e6.7} and \eqref{e6.8};
and
 $$\ba{ll}
(A-\tau I)^\nu W_{n,j}&=n^{-\nu}O\(n|\De s|\)+O(e^{-\de n})+\(\frac{\lam\|A\|}{n+1}\)^{n+1}o(1)\\[1ex]&=O(n^{-\nu}).\ea
$$
Therefore we find that
\be\label{e:7.24}
\phi_n''(\tau)=O(n^{-2(\nu-1)}\,)
\ee
in a small neighborhood $U_n$ of $s$.
\vs We are now in a position  to describe our computation scheme.
\Vs\vs
\noindent {\bf $\bullet$ A combined method computing the  principal eigenvalues}\Vs
Let $A$ be a Perron-like matrix. Take $V=I$ and consider the iteration sequence given by \eqref{e5.6}.
\Vs
{\bf Step 1}. Take a number $N$ suitably large to obtain an approximate  value $s_0$ for $s:=s(A)$  with error estimate $O(N^{-1})$:
$$s_0:=\langle AW_N,\,W_N\rangle.$$ Meanwhile, pick a number  $j\in J$ with $\|W_{N,j}\|=\max_{i\in J}\|W_{N,i}\|.$

\vs
{\bf Step 2}. Determine  the order $\nu$ of $y_j$ via the method in Section \ref{s:6}, which is actually equal to  the cyclic order of $\mb{GE}_s(A)$.

\vs
{\bf Step 3}. Take a number $n<<N$ and solve the  initial value problem
\be\label{e:7.25}
\frac{d\tau}{dt}=-(\gam n)^{2(\nu-1)}\phi_n'(\tau),\hs \tau(0)=s_0
\ee
with parameter $\gam$ to obtain an approximate value $\tau(t_0)$ of the minimum point $\tau_*$ of the function
$\phi_n(\tau)$ given in \eqref{e:7.26}.
 By virtue of Lemma \ref{l:7.1} we have
\be\label{e:7.27}|\tau_*-s|=O\(\|W_{n,j}-q_j(n)\|\).\ee

On the other hand, since $s_0$ is in a small neighborhood of $s$ (hence $s_0$ is close to $\tau_*$),  by \eqref{e:7.15} and the basic knowledge in the qualitative theory of ODEs we deduce that $\tau(t)\ra\tau_*$ exponentially. This may allow us to get a numerical value $\bar s:=\tau(t_0)$ of $\tau_*$ with $|\bar s-\tau_*|<\epsilon/2$  quickly, where $\epsilon$ denotes a tolerance error.  Combining this with  \eqref{e:7.27},
by \eqref{e:7.20} it yields
$$\ba{ll}
|\bar s-s|&=\epsilon/2+O\(\|W_{n,j}-q_j(n)\|\)\\[1ex]
 &\leq \epsilon/2+B_0e^{-\de n}+\(\frac{\lam\|A\|}{n+1}\)^{n+1}o(1).\ea
$$
Since the second and  third terms in the righthand side of the above estimates tend to $0$ exponentially as $n\ra\8$,  it is desirable that a smaller choice of $n$ can produce  a numerical value $\bar s$ of $s$ with desired  accuracy.

\br\label{r:7.2} The role of the term $(\gam n)^{2(\nu-1)}$ in \eqref{e:7.25} is to speed up the convergence of $\tau(t)$ to $\tau_*$, which is due to the simple observation made  in \eqref{e:7.24}. Theoretically, the larger the parameter $\gam$ is taken, the faster the convergence is. But in practice, due to the discretization of the equation and   the existence of round errors,  large $\gam$ may cause serious oscillations  in the computation of $\tau(t)$.
\er

\br\label{r:7.3} The reason why we pick another smaller number $n<<N$ in Step 3 is that the function $\phi_n(\tau)$ may have multiple minimum points, with each one being an attractor of the equation in \eqref{etau}. For large $n$, these minimum points can be very close to each other. Therefore the initial value $s_0$ in \eqref{etau} may fall into the attraction region of a wrong  minimum point. One way to avoid this risk in case the number $n$ in \eqref{etau} is fixed is to enlarge $N$ in Step 1 to obtain a better approximate value $s_0$ so that it is in the attraction region of the right minimum point $\tau_*$.

\er
\subsection{A numerical example}
\vs
\noindent{\bf Example 7.1.} We continue our business  with  the matrix $A$ used in Examples 5.3 and 6.1.
Take  $N=100$. Then from Table \ref{table:6.1} we get that $s_0:=s_{N}=2.0203$. We also infer from the computations involved in Examples 5.3 and 6.1 that   $\|W_{N,4}\|=\max_{i\in J}\|W_{N,i}\|$ (hence $j=4$) and $\nu=3$.

Now the equation \eqref{e:7.25} can be expressed as
\be\label{etau}
\frac{d\tau}{dt}=(\gam n)^{4}\cdot6\cdot\left\langle (A-\tau I)^{2}W_{n,j},\,(A-\tau I)^{3}W_{n,j}\right\rangle,\hs \tau(0)=s_0
\ee
Taking    $\gam=0.2$, $n=20$ and solving \eqref{etau}, we obtain the following table of numerical results, where $\epsilon(\bar s,s)$ denotes the error estimate between $\bar s=\tau(t)$ and the principal eigenvalue  $s$:

\begin{table}[H]\label{table-7.1}
\caption{Numerical results for $\bar s:=\tau(t)$ with $\gam=0.2$, $n=20$}
\begin{center}
\begin{tabular}{|c|>{\centering\arraybackslash}p{3cm}|>{\centering\arraybackslash}p{3cm}|}
\hline
$t$&$\bar s$&$\epsilon(\bar s,s)$\\[1ex]
\hline
0&2.020300000000000& $2.0300\X 10^{-2}$\\
\hline
10&2.013637560931887&$1.3638\X 10^{-2}$\\
\hline
20&2.006802688403743&$6.8027\X 10^{-3}$ \\
\hline
30&2.002493148980474&$2.4931\X 10^{-3}$\\
\hline
40&2.000760076153716&$7.6008\X 10^{-4}$\\
\hline
50&2.000216906077931&$2.1691\X 10^{-4}$\\
\hline
60&2.000060858938304&$6.0859\X 10^{-5}$\\
\hline
70&2.000017020006746&$1.7020\X 10^{-5}$\\
\hline
80&2.000004774174110&$4.7742\X 10^{-6}$\\
\hline
90&2.000001358699798&$1.3587\X 10^{-6}$\\
\hline
100&2.000000406495274&$4.0650\X 10^{-7}$\\
\hline
\end{tabular}
\end{center}
\end{table}

\br Since $s=2$ is an eigenvalue of $A$ with algebraic multiplicity $3$, using the iterative method in Section \ref{s:5}, one can only obtain a numerical value $\bar s$ of $s$ with an error estimate $2.0300\X 10^{-2}$ by taking $N=100$. However, the   above table indicates that if one uses the combined method developed in this section,  the same choice of $N$ can yield  a numerical value $\bar s$ of $s$ with an error estimate $4.0650\X 10^{-7}$.
\er

\br As we have mentioned in Remark \ref{r:7.2},  large values of the parameter $\gam$ in  \eqref{etau} may cause oscillations in the numerical computation of the solution $\tau(t)$ of the equation. This can be seen from the following plots of numerical  solutions $\tau(t)$ of \eqref{etau} corresponding to different choices of $\gam$.
\er

\begin{figure}[H]
\centering
\includegraphics[height=7cm,width=12cm]{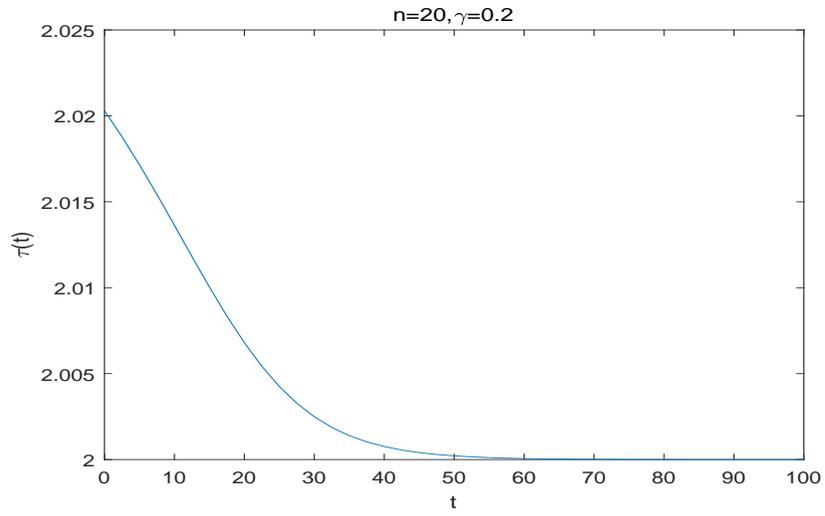}
\caption{Plot of $\tau(t)$ with $n=20$, $\gam=0.2$.}
\label{picture-0.2}
\end{figure}

\begin{figure}[H]
\centering
\includegraphics[height=7cm,width=12cm]{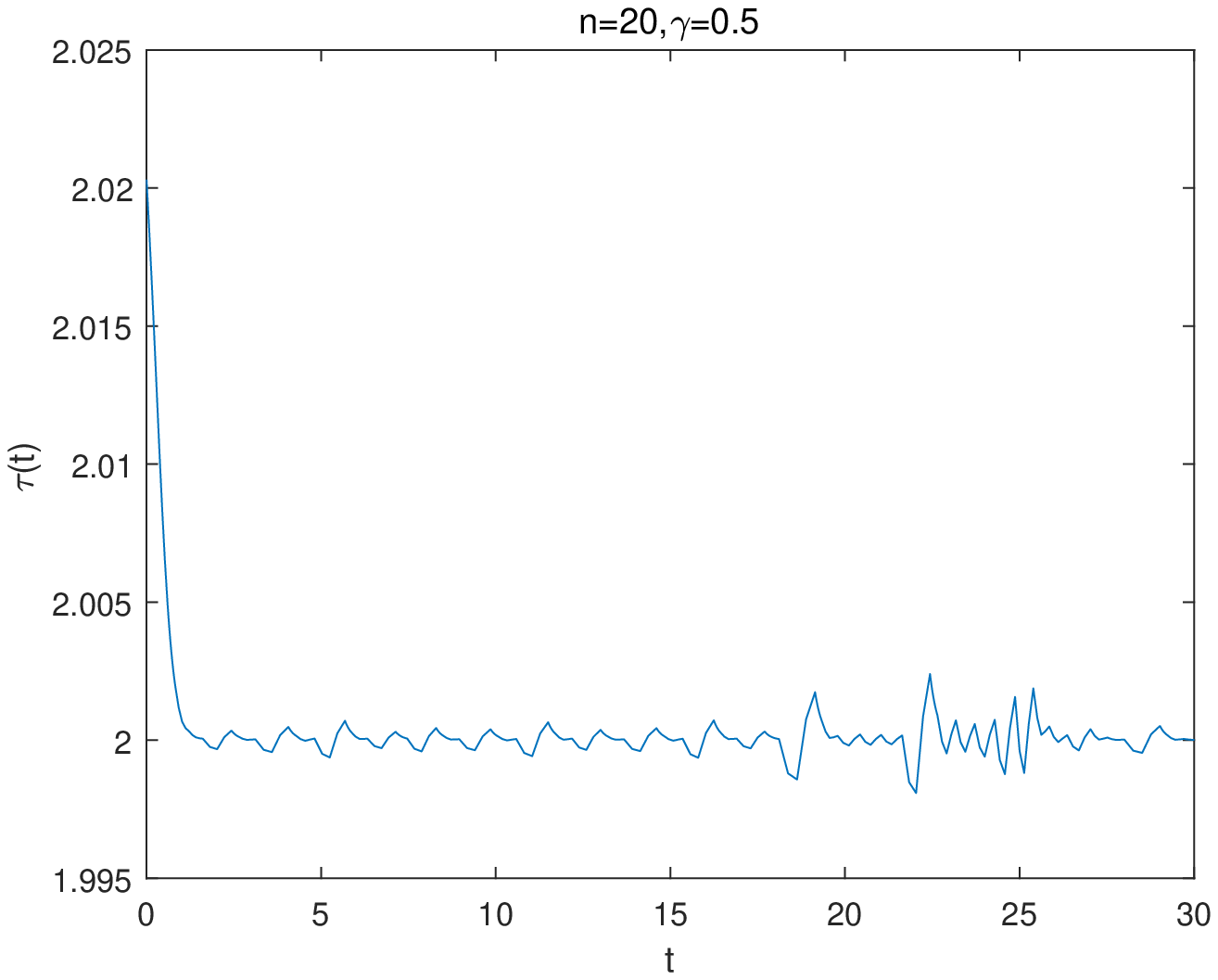}
\caption{Plot of $\tau(t)$ with $n=20$, $\gam=0.5$}
\label{picture-0.5}
\end{figure}

\begin{figure}[H]
\centering
\includegraphics[height=7cm,width=12cm]{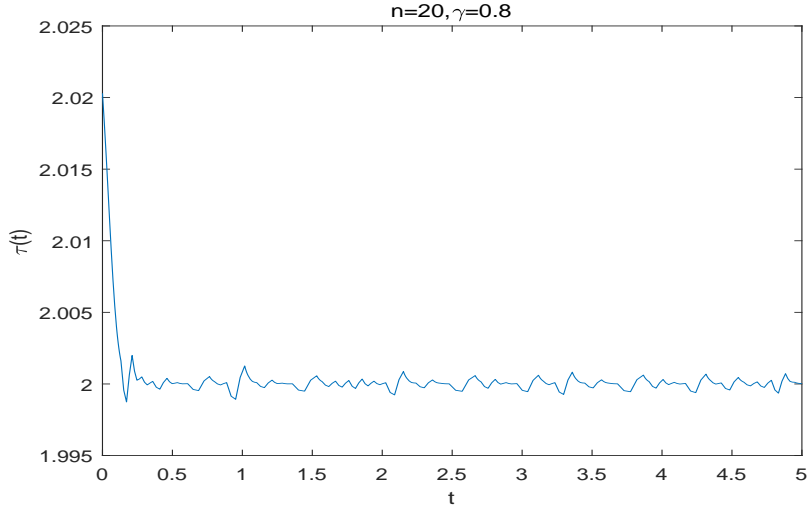}
\caption{Plot of $\tau(t)$ with $n=20$, $\gam=0.8$}
\label{picture-0.8}
\end{figure}

\section{The Computation of  Principal Generalized  Eigenspaces}\label{s:8}
\vs
Let $A\in\bbM_m$ be a Perron-like matrix with principal eigenvalue $s$. Given a nonsingular matrix $V=[v_i]_{i\in J}$, let $M_n(k)$ ($n,k\geq 0$) be the iteration sequence defined in \eqref{e5.6}. If $s$ is semisimple, then by Theorem \ref{t:4.1'} $M_n(n)$ converges to $Y_0:=Y/\|Y\|$ exponentially, where $Y=[y_i]_{i\in J}$ with $y_i=\Pi_1 v_i$, and $\Pi_1$ is the projection from $E:=\R^m$ to $\mb{GE}_s(A)$. Since $\mb{GE}_s(A)=\mb{E}_s(A)=\mb{span}\{y_i:\,\,i\in J\}$ (see Remark \ref{r:3.5}), the computation of the invariant subspace $\mb{GE}_s(A)$ has already  been solved (at least in theoretical sense).  In the  non-semisimple case,  we infer from Theorem \ref{t:4.0} that $W_n:=M_n(n)$ converges to a matrix $\Xi$ whose column vectors can only span the   principal dominant  eigenspace $\mb{DE}_s(A)$ which is  a proper  subspace of $\mb{GE}_s(A)$.

Because  in any case we always have
\be\label{e:6.5}\mb{GE}_s(A)=\mb{span}\{y_i:\,\,i\in J\},\ee  one can select a  well-conditioned basis for $\mb{GE}_s(A)$ from the column vectors  of $Y$. Therefore the computation of $\mb{GE}_s(A)$ can be transformed into that of the matrix $Y$ or related matrices.

\subsection{Some fundamental  convergence results}
\vs
For convenience, denote by $\mb{MGE}_s(A)$ the set of $m\X m$ matrices whose column vectors are in  $\mb{GE}_s(A)$. Since $A\(\mb{GE}_s(A)\)\subset \mb{GE}_s(A)$, we clearly have  \be\label{e:8.2}AM\in \mb{MGE}_s(A),\hs \A\,M\in \mb{MGE}_s(A).\ee  In what follows,  $C$ denotes  a general positive constant which may be different from one to another.


 Write $V=Y+Z$. Then    $Y\in\mb{MGE}_s(A)$,   and $Z=[z_i]_{i\in J}$ is a matrix with column vectors $z_i$ in $ E_2:=\bigoplus_{\mu\in \sig(A)\sm\{s\}}\mb{GE}_\mu(A)$.

\vs Let $
\hat X(t):=e^{t(A-sI)}V$. We observe that
$$\hat X(t)=e^{t(A-sI)}Y+e^{t(A-sI)}Z=Q_{Y}(t)+e^{-st}e^{tA}Z:=Q_{Y}(t)+\hat Z(t),$$
{where } $Q_Y(t)=\[Q_{y_i}(t)\]_{i\in J}$.
 The same calculations in  \eqref{e3.8} apply to show that
\be\label{e:6.16}\|\hat Z(t)\|\leq C\|V\|e^{-\de t},\hs t\geq 0.\ee

By the definition of the  characteristic polynomials of generalized eigenvectors (see \eqref{ecp}) it is easy to  see  that $Q_{y_i}(t)\in \mb{GE}_s(A)$ for all $i\in J$. Consequently
$$(A-s I)^kQ_Y(t)=0,\hs k\ge \nu:=\mb{Cord}\(\mb{GE}_s(A)\).$$
Hence  we deduce that
\be\label{es3.42}\ba{ll}
Y&=e^{-t(A-sI)}e^{t(A-sI)}Y=e^{-t(A-sI)}Q_{Y}(t)\\[2ex]
&=\sum_{k=0}^{\8}(-1)^k\,\frac{t^k}{k!}(A-sI)^kQ_{Y}(t)\\[2ex]
&=\sum_{k=0}^{d}(-1)^k\,\frac{t^k}{k!}(A-sI)^kQ_{Y}(t),\hs t\geq0\ea
\ee
for any integer $d$ with $\nu-1\leq d\leq m-1$. Inspired by this  observation, we take an integer   $d$ with $\nu-1\leq d\leq m-1$ and  define
\be\label{e:dp}P(t)=\sum_{k=0}^d(-1)^k\,\,\frac{t^k}{k!}(A-sI)^k. \ee
Then \eqref{es3.42} can be written as
\be\label{e:6.12}
Y=P(t)Q_Y(t),\hs t\geq 0.
\ee

\bt\label{t:6.1} There is a constant $C>0$ such that
\be\label{e:6.3}\|P(t)\hat X(t)-Y\| \leq C\|V\|(1+t^d) e^{-\de t},\hs t\geq0.
\ee
\et
{\bf Proof.} By  \eqref{e:6.16} and \eqref{e:6.12}  we deduce   that
$$\ba{ll}
\|P(t)\hat X(t)-Y\|&=\|P(t)\hat X(t)-P(t)Q_{Y}(t)\|\\[2ex]
&\leq \|P(t)\|\|\hat Z(t)\|\leq C\|V\|(1+t^d)e^{-\de t},\hs t\geq 0.\ea
$$
This is precisely what we desired. $\bx$

\vs
In practice,   one may  have  only  a numerical  value $\ol s:=s+\De s$ of $s$. Therefore   instead of $\hat X(t)$ and the polynomials $P(t)$ defined in \eqref{e:dp}, we need to take into account   their  approximations:
\be\label{e:8.6}
\bar X(t)=e^{t(A-\bar sI)}V,\hs \bar P(t)=\sum_{k=0}^{d}(-1)^k\,\,\frac{t^k}{k!}(A-\bar sI)^k.
\ee
Note  that
\be\label{e:8.5}\ba{ll}
\bar X(t)&=e^{t(A-\bar s I)}V=e^{t(A-\bar s I)}Y+e^{t(A-\bar s I)}Z\\[2ex]
&=e^{-t\De s}Q_Y(t)+\bar Z(t):=\bar Q_Y(t)+\bar Z(t).\ea
\ee
Hence
\be\label{e:8.8}
\bar P(t) \bar X(t)=\bar P(t)\bar Q_Y(t)+\bar P(t)\bar Z(t)=\bar Y(t)+\bar P(t)\bar Z(t),
\ee
where
\be\label{e:8.9}
\bar Y(t)=\bar P(t)\bar Q_Y(t).
\ee
Since $Q_Y(t)\in \mb{MGE}_s(A)$ and $\bar P(t)$ is a polynomial of $A$,   by \eqref{e:8.2} it is easy to see that
\be\label{e:6.15}
\bar Y(t)=\bar P(t)\bar Q_Y(t)=e^{-t\De s}\bar P(t) Q_Y(t) \in \mb{MGE}_s(A),\hs t\geq 0.\ee

\vs We may assume in advance that  $|\De s|<\de/3$. Then by \eqref{e:6.16} one has
\be\label{e:8.22}
\|\bar Z(t)\|=\|e^{t(A-\bar s I)}Z\|\leq e^{-t\De s}\|e^{t(A-s I)}Z\|\leq C\|V\|e^{-\frac{2}{3}\de t},\hs t\geq 0.
\ee
Thus one deduces  that
\be\label{e:8.10}
\|\bar P(t)\bar Z(t)\|\leq \|\bar P(t)\|\|\bar Z(t)\|\leq C\|V\|(1+t^d)e^{-\frac{2}{3}\de t},\hs t\geq 0.
\ee
It follows by \eqref{e:8.8} that
$$
\|\bar P(t)\bar X(t)-\bar Y(t)\|\leq  C\|V\|(1+t^d)e^{-\frac{2}{3}\de t},\hs t\geq 0.
$$
\vs
We also observe that
\be\label{e:8.11}\ba{ll}
\|\bar Y(t)-Y\|&=(\mb{by \eqref{e:6.12}})=\|\bar P(t)\bar Q_Y(t)-P(t)Q_Y(t)\|\\[2ex]
&=(\mb{by \eqref{e:8.5}})=\left\|e^{-t\De s}\bar P(t) Q_Y(t)-P(t)Q_Y(t)\right\|\\[2ex]
&\leq \|e^{-t\De s}\bar P(t) -P(t)\|\|Q_Y(t)\|.
\ea
\ee
Now we assume  that $|t\De s|\leq 1$. Then
$$\ba{ll}
 \|e^{-t\De s}\bar P(t) -P(t)\|&\leq  \|\(e^{-t\De s}-1\)\bar P(t)\|+ \|\bar P(t) -P(t)\|\\[2ex]
 &\leq C(1+t^d)|t\De s|+C(1+t^d)|\De s|\\[2ex]
 &\leq C(1+t^{d+1})\,|\De s|,\hs t\geq 0.
 \ea
 $$
 Hence by \eqref{e:8.11} we conclude that
$$
\|\bar Y(t)-Y\|\leq C(1+t^{d+1})\,|\De s|\|Q_Y(t)\|\leq C\|V\|(1+t^{d+\nu})\,|\De s|.
$$

Let us write
$$\widetilde{Y}(t):=\bar P(t)\bar X(t)$$ and summarize the above results in the following theorem:

 \bt\label{t:8.2} Suppose  $|\De s|<\de/3$ and that $t|\De s|\leq 1$. Then there is a constant $C>0$ depending only upon $\|A\|$ and $d$ such that
 \be\label{e:8.12}
 \|\widetilde{Y}(t)-\bar Y(t)\|\leq  C\|V\|(1+t^d)e^{-\frac{2}{3}\de t},\hs t\geq 0,
 \ee  and
\be\label{e:8.13}\|\bar Y(t)-Y\| \leq C\|V\|(1+t^{d+\nu})\,|\De s|,\hs t\geq0.
\ee
\et

\vs
We infer from  \eqref{e:6.5}  that there are $i_1,i_2,\cdots,i_\ell\in J$ such that $\{y_{i_k}\}_{1\leq k\leq\ell}$ forms a basis of $\mb{GE}_s(A)$. Hence
$$
\|a_1y_{i_1}+\cdots+a_\ell y_{i_{\ell}}\|\geq \eta>0,\hs\A a=(a_1,\cdots,a_\ell)^{\mb{\tiny T}}\in \R^\ell,\,\,\|a\|=1.
$$
By \eqref{e:8.13} we see that for each fixed pair $(t,\De s)$, if $\De s$ is sufficiently small so that $\|\bar Y (t)-Y\| \leq \eta/2$, then one has
$$
\|a_1\bar y_{i_1}(t)+\cdots+a_\ell \bar y_{i_{\ell}}(t)\|\geq \eta/2,\hs\A a=(a_1,\cdots,a_\ell)^{\mb{\tiny T}}\in \R^\ell,\,\,\|a\|=1,
$$
where $\bar y_i(t)$ denotes the $i$-th column vector of the matrix $\bar Y (t)$. This indicates that the vectors $\bar y_{i_k}(t)\,\,({1\leq k\leq \ell})$ are fairly linearly independent and therefore, in view of \eqref{e:6.15}, form a well-conditioned basis of the space $\mb{GE}_s(A)$.

Unfortunately in practice  $\bar Y (t)$ is generally  unknown. However,  by Theorem \ref{t:8.2} there exists $t_0>0$ such that \be\label{e:8.15}\|\widetilde{Y}(t)-\bar Y(t)\|<\epsilon\ee  for $t\geq t_0$ and $|\De s|<\de/3$, where $\epsilon$ denotes tolerance error.
Suppose now that $\De s$ is sufficiently small so that we can pick  $t\geq t_0$ suitably  large such that $\|\bar Y(t)-Y\|<<1$.  Then by \eqref{e:8.15} one can select a well-conditioned numerical basis  for $\mb{GE}_s(A)$ from the column vectors of $\widetilde{Y}(t)$.

\subsection{An iterative method for the computation of $\bar{Y}(t)$} \vs
The above  discussion reduces the computation of  the space $\mb{GE}_s(A)$ to that of the matrix $\bar{Y}(t)$ for $t$ sufficiently large and $|\De s|$ sufficiently small, which can be done by simply computing $\widetilde{Y}(t)=\bar P(t)\bar X(t)$. Note that $\bar X(t)$ is a matrix exponential.
Although there have been many excellent computation methods for matrix exponentials, here we prefer to develop an iterative scheme as the one in Section \ref{s:5} to compute the normalization of $\widetilde{Y}(t)$, which avoid some technical difficulties (such as overflow for large $t$) in the computation.

\Vs
\noindent $\bullet$ {\bf Some basic estimates} \,\,Let us  first give some estimates on the lower bounds of the norms of $\bar Q_Y(t)$ and $\bar X(t)$.
Since $$Q_Y(t)=\sum_{k=0}^{\nu-1}\frac{t^k}{k!}(A-s I)^k Y=e^{t(A-s I)}Y,$$
we deduce that $\|Q_Y(t)\|>0$ for all $t\geq 0$. On the other hand, by the first equality in the above equation it is easy to see that
$$
\|Q_Y(t)\|=t^{\nu-1}\(\frac{1}{(\nu-1)!}\|(A-s I)^{\nu-1}Y\|+O(t^{-1})\)\hs\mb{as }t\ra\8.
$$
Thus there is a $t_0>1$ such that
$$
\|Q_Y(t)\|\geq c_0t^{\nu-1},\hs t>t_0
$$
for some $c_0>0$. Take a $c_1>0$ with $c_1<c_0$ such that $\|Q_Y(t)\|\geq c_1$ for $t\in[0,t_0]$. Then
\be\label{e:8.20}
\|Q_Y(t)\|\geq c_1>0,\hs \A\,t\geq 0.
\ee
Consequently
\be\label{e:8.21}
\|\bar Q_Y(t)\|=e^{-t\De s}\|Q_Y(t)\|\geq c_1e^{-t|\De s|},\hs \A\,t\geq 0.
\ee

Further by \eqref{e:8.5} and \eqref{e:8.22} we conclude that
\be\label{e:8.23}
\|\bar X(t)\|\geq \|\bar Q_Y(t)\|-\|\bar Z(t)\|\geq c_1e^{-t|\De s|}-C\|V\|e^{-\frac{2}{3}\de t}
\ee
for $t\geq 0$. Take a $t_1>1$ such that $C\|V\|e^{-\frac{\de}{3}t}<c_1/2$ for $t>t_1$. Then since $|\De s|<\de/3$, we have
$$
\|\bar X(t)\|\geq c_1e^{-t|\De s|}-\frac{1}{2}c_1e^{-\frac{\de}{3} t}\geq \frac{1}{2}c_1e^{-t|\De s|},\hs t>t_1.
$$
Noticing that $\bar X(t)=e^{t(A-\bar s I)}V\ne 0$ for all $t\geq 0$, as in the case of $Q_Y(t)$, we can pick a number
$c_2>0$ with $c_2<c_1/2$ such that $$\|\bar X(t)\|\geq c_2e^{-t|\De s|}$$ for $t\in[0,t_1]$. Then
\be\label{e:8.24}
\|\bar X(t)\|\geq  c_2e^{-t|\De s|},\hs t\geq 0.\ee

\vs
Now let us evaluate  $\left\|\bar P\(\frac{\bar X}{\|\bar X\|}\)-\frac{\bar Y}{\|\bar Q_Y\|}\right\|$, where
$$\bar P=\bar P(t),\hs \bar X=\bar X(t),\hs \bar Y=\bar Y(t),\hs \bar Q_Y=\bar Q_Y(t).$$
Recalling that $\bar X=\bar Q_Y+\bar Z$ and $\bar Y=\bar P \bar Q_Y$, we have
\be\label{e:8.25}\ba{ll}
\left\|\bar P\(\frac{\bar X}{\|\bar X\|}\)-\frac{\bar Y}{\|\bar Q_Y\|}\right\|&=\left\|\frac{\bar Y}{\|\bar X\|}+\frac{\bar P\bar Z}{\|\bar X\|}-\frac{\bar Y}{\|\bar Q_Y\|}\right\|\\[2ex]
&\leq \left\|\frac{\bar Y}{\|\bar X\|}-\frac{\bar Y}{\|\bar Q_Y\|}\right\|+\frac{\|\bar P\bar Z\|}{\|\bar X\|}
\ea
\ee

By \eqref{e:8.22} and \eqref{e:8.24} we easily deduce that
\be\label{e:8.26}
\frac{\|\bar P\bar Z\|}{\|\bar X\|}\leq C\|V\|(1+t^d)e^{-\frac{\de}{3}t},\hs t\geq 0.
\ee

If we write $\|\bar X\|=\|\bar Q_Y\|+r(t)$, then by $\bar X=\bar Q_Y+\bar Z$ we see that
\be\label{e:8.27}
r(t)\leq \|\bar X-\bar Q_Y\|=\|\bar Z\|\leq C\|V\|e^{-\frac{2}{3}\de t}.
\ee
Thus
$$
\left\|\frac{\bar Y}{\|\bar X\|}-\frac{\bar Y}{\|\bar Q_Y\|}\right\|\leq \left|\frac{1}{\|\bar X\|}-\frac{1}{\|\bar Q_Y\|}\right|\|\bar Y\|\leq \frac{r(t)}{\|\bar X\|\|\bar Q_Y\|}\|\bar Y\|
$$
We observe that
$$
\|\bar Y\|= \|\bar P\bar Q_Y\|\leq \|\bar P\|\|e^{-t\De s} Q_Y\|\leq C\|Y\|(1+t^{d+\nu-1})e^{t|\De s|}
$$
for $t\geq 0$. Hence by \eqref{e:8.27} and the estimates on the lower bounds of $\|\bar X\|$ and $\|\bar Q_Y\|$ we deduce that
$$
\left\|\frac{\bar Y}{\|\bar X\|}-\frac{\bar Y}{\|\bar Q_Y\|}\right\|\leq C\|V\|\|Y\|(1+t^{d+\nu-1})e^{-\frac{\de}{3}t},\hs t\geq 0
$$
provided that $|\De s|<\de/9$. Combining this with \eqref{e:8.25} and \eqref{e:8.26}, we have
\bl\label{l:8.3} Assume that $|\De s|<\de/9$. Then there exist $C=C_V>0$ such that
\be\label{e:8.28}\ba{ll}
\left\|\bar P\(\frac{\bar X}{\|\bar X\|}\)-\frac{\bar Y}{\|\bar Q_Y\|}\right\|\leq C(1+t^{d+\nu-1})e^{-\frac{\de}{3}t},\hs t\geq 0.
\ea
\ee
\el
\Vs

\noindent$\bullet$ {\bf The iterative method} \,\,Let $\bar A:=A-\bar s I$.
For each fixed $n\in \bbN$,  set $$T_n=\sum_{k=0}^n\frac{1}{k!}\bar A^k.
$$
Since $T:=e^{t\bar A}$ is nonsingular, we deduce  that $T_n$ is nonsingular as well  provided that  $n$ is sufficiently large.
Define a mapping  $K_n$ on $$\Sigma_1:=\{M\in \bbM_m:\,\,\|M\|=1\}$$ as follows:
$$K_n M=\frac{T_n M}{\|T_n M\| },\hs \A\,M\in \Sigma_1.$$
Given   a nonsingular matrix   $V=[v_i]_{i\in J}\in \Sigma_1$, define an iteration sequence  as
\be\label{e:8.29}
\cS_n(0)=V,\,\,\,\,\cS_n(k+1)=K_n \cS_n(k),\hs\,\, k\in\bbN.
\ee
Set $S_n:=\cS_n(n).$
\bt\label{t:8.5}Suppose that $|\De s|<\de/9$.   Then for $ n\in \bbN$, we have
\be\label{e:8.35}\ba{ll}
\left\|\frac{\bar P(n) S_n}{\|\bar P(n) S_n\|}-\frac{\bar Y(n)}{\|\bar Y(n)\|}\right\|&\leq C\|\bar P(n)S_n\|^{-1}(1+n^d)\(\frac{e \bar\lam\|\bar A\|}{n+1}\)^{n+1}o(1)+\\[2ex]
&\hs +C\|\bar Y(n)\|^{-1}(1+n^{d+2(\nu-1)})e^{-\frac{2}{9}n\de}.
\ea
\ee
\et
{\bf Proof.} Let $\bar S_n={\bar X(n)}/{\|\bar X(n)\|}$. Repeating the same argument in Section \ref{s:5.1} with the matrix $A$ therein replaced by $\bar A$, it can be shown that
\be\label{e:8.31}
\|S_n-\bar S_n\|\leq \bar \lam^{n+1} o_n,
\ee
where $o_n=\frac{\|\bar A\|^{n+1}}{(n+1)!}$. Therefore
\be\label{e:8.32}
\|\bar P(n)S_n-\bar P(n)\bar S_n\|\leq C(1+n^d)\bar \lam^{n+1} o_n.
\ee

For convenience in statement, in what follows we assign $\frac{O_{m\X m}}{\|O_{m\X m}\|}=I$ for the zero matrix $O_{m\X m}\in\bbM_m$, and  assign $\frac{0}{0}=1$ for the number  $0\in\R$. Then one trivially verifies that
$$
\left\|\frac{M_1}{\|M_1\|}-\frac{M_2}{\|M_2\|}\right\|\leq \frac{2}{\|M_2\|}\|M_1-M_2\|,\hs \A\,M_1,M_2\in \bbM.$$
 Let $F_n=\bar Y(n)/\|\bar Q_Y(n)\|$.
Since $\|\bar Q_Y(t)\|\leq C(1+t^{\nu-1})e^{t|\De s|}$ for $t\geq 0$, we have
$$
\left\|F_n\right\|\geq \frac{\|\bar Y(n)\|}{C(1+n^{\nu-1})e^{t|\De s|}},\hs n\in \bbN.
$$
Therefore
\be\label{e:8.33}\ba{ll}
\left\|\frac{\bar P(n)\bar S_n}{\|\bar P(n)\bar S_n\|}-\frac{\bar Y(n)}{\|\bar Y(n)\|}\right\|&=\left\|\frac{\bar P(n)\bar S_n}{\|\bar P(n)\bar S_n\|}-\frac{\bar F_n}{\|\bar F_n\|}\right\|\\[2ex]
&\leq 2\|F_n\|^{-1}\left\|{\bar P(n)\bar S_n}-{\bar F_n}\right\|\\[2ex]
&\leq C\|\bar Y(n)\|^{-1}(1+n^{d+2(\nu-1)})e^{-\frac{2}{9}n\de}
\ea
\ee
for all $n\in \bbN$.

Similarly we have
\be\label{e:8.34}\ba{ll}
\left\|\frac{\bar P(n) S_n}{\|\bar P(n) S_n\|}-\frac{\bar P(n)\bar S_n}{\|\bar P(n)\bar S_n\|}\right\|&\leq C\|\bar P(n)S_n\|^{-1}\left\|{\bar P(n) S_n}-{\bar P(n)\bar S_n}\right\|\\[2ex]
&\leq (\mb{by }\eqref{e:8.32})\\[2ex]
&\leq  C\|\bar P(n) S_n\|^{-1}(1+n^d){\bar \lam}^{n+1}o_n\\[2ex]
&=C\|\bar P(n) S_n\|^{-1}(1+n^d)\(\frac{e \bar\lam\|\bar A\|}{n+1}\)^{n+1}o(1)
\ea
\ee
for $n\in\bbN$. Here we have used the Stirling's formula. Combining this with \eqref{e:8.33} one immediately obtains the validity of \eqref{e:8.35}. $\bx$

\br\label{r:8.3} In many cases one can use the method in  Section \ref{s:6} to compute   the cyclic order $\nu:=\mb{\em Cord}\(\mb{\em GE}_s(A)\)$. Whence  $\nu$ is determined, the best choice of the number $d$ is naturally  $d=\nu-1$. 
\er

\br As we have pointed out in Remark \ref{r:5.5}, in practice one may  take an arbitrary nonsingular matrix $V$ as the initial one in  the iteration.
\er
\Vs

\noindent{$\bullet$} {\bf Computation scheme of }\,$B_n:=\bar{Y}(n)/\|\bar{Y}(n)\|$
\vs
As we have seen, the set of column vectors of $B_n$ provides  a well-conditioned spanning set for the principal generalized eigenspace $\mb{GE}_s(A)$. Thanks to Theorem \ref{t:8.5}, it is desirable to get numerical values of $B_n$ within tolerance error quickly by using the iterative method developed here. The computation scheme is as follows.
\Vs
{\bf Step 1}. Take $V=I$. Use the method in Section \ref{s:6} to determine the  cyclic order $\nu=\mb{Cord}\(\mb{GE}_s(A)\)$.
\vs
{\bf Step 2}. Use the combined method in Section \ref{s:7} to derive a  numerical value $\bar s$ of $s$ as accurate as possible.
\vs
{\bf Step 3}. Use the iteration sequence defined in \eqref{e:8.29} to compute $$\widetilde{B}_n:={\bar P(n) S_n}\,/{\|\bar P(n) S_n\|},$$
where $S_n=\cS_n(n)$. In general $\widetilde{B}_n$ approaches $B_n$ exponentially. Hence the desired accuracy may  be obtained swiftly by taking a suitably large $n$.

\Vs

\br The invariant subspaces of a matrix can be successively computed via  the famous Krylov subspace methods and its combinations with other algorithms; see e.g. Watkins \cite{Wat2} for a systematic discussion on this topic. In Bai and Demmel \cite{BaiZJ} the authors  developed another method to compute invariant subspaces by using matrix sign functions.
\er

\subsection{A numerical example}
\vs
\noindent{\bf Example} 8.1. Let
$$A=\(\begin{matrix}-1&\frac{1}{2}&0&3&0&-\frac{1}{2}&3\\
                   0&1&2&0&-2&1&0\\
                   3&\frac{1}{2}&-1&-3&3&-\frac{1}{2}&-3\\
                   1&1&-2&1&2&-1&-1\\
                   0&\frac{1}{2}&-2&0&4&-\frac{1}{2}&0\\
                   0&-1&1&0&-2&3&-1\\
                   -3&-\frac{1}{2}&2&3&-2&\frac{1}{2}&5\end{matrix}\)$$
                   The principal  eigenvalue of the matrix is  $s=2$ with an  algebraic multiplicity $5$ and a geometric multiplicity $3$. Taking $V=I=[e_i]_{i\in J}$, where $e_i$ is the $i$-th column vector of $I$, $J=\{1,2,\cdots,7\}$, we have
                  $$Y=\[y_i\]_{i\in J}=\(\begin{matrix}-\frac{13}{2}&0&-\frac{1}{2}&\frac{15}{2}&\frac{1}{2}&0&\frac{15}{2}\\
                   15&1&3&-15&-3&0&-15\\
                   -\frac{21}{2}&0&-\frac{5}{2}&\frac{21}{2}&\frac{7}{2}&0&\frac{21}{2}\\
                   -14&0&-3&15&3&0&14\\
                   -\frac{21}{2}&0&-\frac{5}{2}&\frac{21}{2}&\frac{7}{2}&0&\frac{21}{2}\\
                   9&0&2&-9&-2&1&-9\\
                   \frac{15}{2}&0&\frac{5}{2}&-\frac{15}{2}&-\frac{5}{2}&0&-\frac{13}{2}\end{matrix}\)$$
 where $y_i=\Pi_1 e_i$, and $\Pi_1$ is the projection from $\R^7$ to $\mb{GE}_s(A)$. It is known that
$$
\mb{GE}_s(A)=\mb{span}\{y_i:\,\,i\in J\}.
$$
\Vs
\noindent$\bullet$ {\bf Numerical simulation}
\Vs
                   (1) Taking $N=100$ and using  the method in Section \ref{s:5}, we can get an approximate value $s_N$ of $s$: $s_N=2.0203$.
                    \vs
                    (2)  \,We use the method in Section \ref{s:6} to determine the cyclic order $\nu:=\mb{Cord}\(\mb{GE}_s(A)\)$ of the space $\mb{GE}_s(A)$.

                     Set $\ve=0.10$, and take $N=100$ and $s_N=2.0203$. Then the number $j$ satisfying \eqref{ej} is $j=7$. Take different $n<<N$ and compute $\b_k(n,s_N)$ defined in \eqref{e:6.48} for $k\in J$, we get the following table:
                    \begin{table}[H]
\caption{Values of $\b_k(n,s_N)$ with $N=100,s_N=2.0203$}
\begin{center}
\begin{tabular}{|c|>{\centering\arraybackslash}p{1.25cm}|>{\centering\arraybackslash}p{1.25cm}|>
{\centering\arraybackslash}p{1.25cm}|>{\centering\arraybackslash}p{1.25cm}|>{\centering\arraybackslash}p{1.25cm}|>
{\centering\arraybackslash}p{1.25cm}|>{\centering\arraybackslash}p{1.25cm}|}
\hline
$n$&$\b_1(n,s_N)$&$\b_2(n,s_N)$&$\b_3(n,s_N)$&$\b_4(n,s_N)$&$\b_5(n,s_N)$&$\b_6(n,s_N)$&$\b_7(n,s_N)$\\[2ex]
\hline
4&3.9192&1.2926&0.2612&14.1965&2.5930&27.7718&87.5299\\
\hline
5&3.8535&1.2355&0.0228&101.9814&10.4822&4.2120&31.0484\\
\hline
6&3.8418&1.1058&0.0475&22.2124&23.1030&14.1662&6.2880\\
\hline
7&3.8398&0.9839&0.1460&2.4422&46.0797&26.8717&19.3435\\
\hline
8&3.8211&0.8760&0.2686&0.2745&117.1914&42.0806&34.7450\\
\hline
\end{tabular}
\end{center}
\end{table}
\noindent It can be seen that if $n=5,\,6$ then the dichotomy property in  \eqref{e:6.54} is fulfilled with $k_0=3$.
Hence we deduce that  $\nu=3$.
\vs
(3) Now we use the combined method in Section \ref{s:7} to obtain a refined approximate value $\bar s$ of $s$. Let $s_0=s_N=2.0203$ and  $j=7$. Consider the initial value problem:
$$
\frac{d\tau}{dt}=(\gam n)^{4}\cdot6\cdot\left\langle(A-\tau I)^{2}W_{n,j},\,(A-\tau I)^{3}W_{n,j}\right\rangle,\hs \tau(0)=s_0,
$$
where $W_{n,j}$ is the $j$-th column vector of $W_n:=M_n(n)$, and $\{M_n(k)\}_{n,k\in\bbN}$ is the sequence  given by \eqref{e5.6}. Taking $\gam=0.20$ and $n=20$, we obtain that
$$\bar{s}=\tau(150)=2.000005037291918.$$

\vs(4) Use the iteration sequence defined in \eqref{e:8.29} to compute
$$\widetilde{B}_n:={\bar P(n) S_n}\,/{\|\bar P(n) S_n\|},\hs \mb{where }S_n=\cS_n(n).$$
Let $B_n:=\bar{Y}(n)/\|\bar{Y}(n)\|$. Denote by $\epsilon(\widetilde{B}_n,\,B_n)$ the error estimate between $\widetilde{B}_n$ and $B_n$, and by  $\epsilon\(B_n,\,Y/\|Y\|\)$ the error estimate between $B_n$ and $Y/\|Y\|$. Then we have the following table.

\begin{table}[H]
\caption{Numerical results}
\begin{center}
\begin{tabular}{|c|>{\centering\arraybackslash}p{2.5cm}|>{\centering\arraybackslash}p{2.5cm}|}
\hline
$n$&$\epsilon(\widetilde{B}_{n},B_{n})$&$\epsilon(B_{n},Y/\|Y\|)$\\[2ex]
\hline
20&$1.9181\times10^{-6}$&0.0015\\
\hline
22&$3.4588\times10^{-7}$&0.0020\\
\hline
25&$2.5222\times10^{-8}$&0.0029\\
\hline
28&$1.5701\times10^{-9}$&0.0041\\
\hline
30&$2.8581\times10^{-10}$&0.0051\\
\hline
\end{tabular}
\end{center}
\end{table}

\section*{References}
\small{

}

\end{document}